\documentclass[amssymb,11pt]{article}
\usepackage{amsthm}
\usepackage{amsmath}
\usepackage{amssymb}
\usepackage{graphicx}
\usepackage{color}
\usepackage{multirow}
\usepackage{longtable}
\usepackage{subcaption}
\usepackage{comment}
\usepackage{url}
\textwidth=16.8cm \textheight=23.5cm \vspace{1.4cm}

\newtheorem{lemma}{Lemma}[section]
\newtheorem{definition}{Definition}[section]

\newtheorem{algorithm}{Algorithm}[section]
\newtheorem{theorem}{Theorem}[section]
\newtheorem{corollary}{Corollary}[section]
\newtheorem{proposition}{Proposition}[section]
\newtheorem{remark}{Remark}[section]
\newtheorem{example}{Example}[section]

\numberwithin{equation}{section}
\newcommand{\beq}{\begin{equation}}
\newcommand{\eeq}{\end{equation}}
\newcommand{\beqa}{\begin{eqnarray}}
\newcommand{\eeqa}{\end{eqnarray}}
\newcommand{\beqas}{\begin{eqnarray*}}
\newcommand{\eeqas}{\end{eqnarray*}}
\newcommand{\ba}{\begin{array}}
\newcommand{\ea}{\end{array}}
\newcommand{\bi}{\begin{itemize}}
\newcommand{\ei}{\end{itemize}}

\newcommand{\nn}{\nonumber}
\DeclareMathOperator*{\Argmin}{Argmin}  
\DeclareMathOperator*{\argmin}{argmin}

\def\bj{{\mathbb{j}}}
\def\bx{{\bar x}}
\def\bd{{\bar d}}
\def\cI{{\cal I}}
\def\cJ{{\cal J}}
\def\cK{{\cal K}}

\def\cN{{\cal N}}
\def\conv{{\rm conv}}
\def\cT{{\cal T}}
\def\hh{{\hat h}}
\def\hphi{{\hat \phi}}
\def\o{{\rm o}}
\def\tF{{\widetilde F}}
\def\tQ{{\widetilde Q}}

\begin{document}
\title{Penalty and Augmented Lagrangian Methods for \\Constrained DC Programming\footnote{This paper has been accepted for publication in Mathematics of Operations Research.}
}
\author{
Zhaosong Lu%
\thanks{
Department of Industrial and Systems Engineering, University of Minnesota, USA (email: {\tt zhaosong@umn.edu}).}
\and
Zhe Sun%
\thanks{
 School of Mathematics and Statistics, Jiangxi Normal University, Nanchang, China 330022. (email: {\tt snzma@126.com}). This author was supported by National Natural Science Foundation of China (Grant No. 11761037 and 11501265) and Natural Science Foundation of Jiangxi Province (Grant No. 20181BAB201009).
}
\and
Zirui Zhou%
\thanks{Huawei Technologies Canada, Burnaby, BC, Canada. (email: {\tt zirui.zhou@huawei.com}).
}
}

\date{January 27, 2020 (Revised: February 24, 2021)}

\maketitle

\begin{abstract}

In this paper we consider a class of structured nonsmooth difference-of-convex (DC) constrained DC program in which the first convex component of the objective and constraints  is the sum of a smooth and nonsmooth functions while their second convex component is the supremum of finitely many convex smooth functions. The existing methods for this problem usually have a weak convergence guarantee or require a feasible initial point. Inspired by the recent work (Math Oper.\ Res.\ 42(1):95--118, 2017 by Pang et al.), in this paper we propose two infeasible methods with strong convergence guarantee for the considered problem. {The first one is a penalty method that consists of finding an approximate D-stationary point of a sequence of penalty subproblems. We show that any feasible accumulation point of the solution sequence generated by such a penalty method is a B-stationary point of the problem under a weakest possible assumption that it satisfies a pointwise Slater constraint qualification (PSCQ). The second one is an augmented Lagrangian (AL) method that consists of finding an approximate D-stationary point of a sequence of AL subproblems. Under the same PSCQ condition as for the penalty method, we show that any feasible accumulation point of the solution sequence generated by such an AL method is a B-stationary point of the problem,  and moreover, it satisfies a KKT type of optimality condition for the problem, together with any accumulation point of the sequence of a set of auxiliary Lagrangian multipliers.
We also propose an efficient successive convex approximation method for computing an approximate D-stationary point of the penalty and AL subproblems.
} Finally, some numerical experiments are conducted to demonstrate the efficiency of our proposed methods.


\bigskip

\end{abstract}

\setcounter{equation}{0}

\section{Introduction}

Difference-of-convex (DC) programs are a class of important optimization problems, which generally minimize an objective
function that is the difference of two convex functions subject to constraints defined by the same type of functions. They
have been studied for several decades in the literature (e.g., see \cite{Hong,Lu,Hong16,Pang,Wen18,Lu18,Liu18,Thi18} and references therein). In this paper we are interested in a DC program in the form of
\begin{equation}\label{e1}
\begin{array}{ll}
\min\limits_{x\in X}& F(x)= \phi_0(x)+\zeta_0(x)-\psi_0(x)\\
\mbox{s.t.} &  \phi_{i}(x)+\zeta_i(x)-\psi_{i}(x)\leq 0, \quad\forall i=1,\ldots,I,
\end{array}
\end{equation}
where
\begin{equation}\label{e2}
\psi_{i}(x)=\max\limits_{1\leq j\leq J_i}\{\psi_{i,j}(x)\}, \quad \forall i=0,1,\ldots, I
\end{equation}
for some integers $J_i$'s, $X\subseteq \Re^n$ is a closed convex set, $\zeta_i$'s are convex and continuous  on an open set $\cal S$ containing $X$, $\phi_i$'s, $\psi_{i,j}$'s are convex and continuously differentiable on
$\cal S$, and moreover,  $\nabla\phi_i$ is Lipschitz continuous with constant $L_i \ge 0$, that is,
\[
\|\nabla\phi_i(x)-\nabla\phi_i(y)\|\leq L_i\|x-y\|, \quad\forall x, y\in X
\]
for all $i=0,1,\ldots, I$. In addition, for convenience we assume throughout this paper that $L_0>0$.\footnote{This assumption is very mild. Indeed, if $L_0=0$, one can replace $\phi_0(x)$ and $\psi_0(x)$ by $\phi_0(x)+\|x\|^2/2$ and $\psi_0(x)+\|x\|^2/2$, respectively. Then the resulting problem is equivalent to the original one but with $L_0=1>0$.}

DC program \eqref{e1} has found numerous applications in signal processing, communications and networks, statistical variable selections, finance, and etc (e.g., see \cite{GaRaCa09,Hong,Alvarado,Hong16,Pang,Ahn,Dong,Thi18,Liu,Cui}).  Also, it has been shown in \cite{GaRaCa09,Li,Ahn,Thi18}
that some widely used sparse optimization models can be equivalently reformulated and solved as \eqref{e1}.  Several methods have been proposed in the literature for solving general DC programs (see \cite{Hong, Lu,Alvarado,Thi14a,Hong16,Pang}). Nevertheless, they face some issues when applied to solve problem \eqref{e1} as mentioned below.

When a feasible point $x^0$ of \eqref{e1} is known, the sequential convex
programming (SCP) method (e.g., see \cite{Hong, Lu,Alvarado,Thi14a,Hong16}) can be applied to \eqref{e1} and it generates iterates $\{x^k\}$ as
follows:
 \begin{equation}\label{convex-subprob1}
\begin{array}{rl}
x^{k+1} \in \Argmin\limits_{x\in X}& \phi_0(x)+\zeta_0(x)-(s^k_{\psi_0})^T x\\
\mbox{s.t.} &  \phi_{i}(x)+\zeta_i(x)-\psi_i(x^k)-(s^k_{\psi_i})^T(x-x^k) \leq 0, \quad\forall i=1,\ldots,I,
\end{array}
\end{equation}
where $s^k_{\psi_i} \in \partial \psi_i(x^k)$ for $k \ge 0$ and $i=0,1,\ldots, I$. Under some suitable constraint qualification, every accumulation point $x^\infty$ of the sequence $\{x^k\}$ is a Karush-Kuhn-Tucker (KKT)  point of the problem, that is, there exists some $\lambda \in \Re^I$ together with $x^\infty$ satisfying the KKT conditions
\beq \label{KKT}
\ba{l}
0\in \nabla \phi_0(x^\infty) + \partial \zeta_0(x^\infty) - \partial \psi_0(x^\infty) +
\sum\limits_{i=1}^I  \lambda_i[\nabla \phi_i(x^\infty) + \partial \zeta_i(x^\infty) - \partial \psi_i(x^\infty)] +\cN_X(x^\infty), \\
\lambda_i \ge 0, \quad \phi_{i}(x^\infty)+\zeta_i(x^\infty)-\psi_{i}(x^\infty)\leq0,\quad \lambda_i  [\phi_{i}(x^\infty)+\zeta_i(x^\infty)-\psi_{i}(x^\infty)] =0,\quad \forall i=1,\ldots,I.
\ea
\eeq
 Though SCP looks quite natural, it encounters some practical issues.  Firstly, the initial feasible point $x^0$ is usually unknown. Secondly, it is typically difficult to find an {\it exact} solution of subproblem \eqref{convex-subprob1}.

 Le Thi et al.\ \cite{Thi14a} proposed two penalty type of methods { (named as DCA1 and DCA2)
 for solving {\it a special case} of problem \eqref{e1} with $\psi_i$ being continuously differentiable on $X$,\footnote{The differentiability of $\psi_i$'s is required for the convergence of the methods (see \cite[Assumption 2]{Thi14a}).} that is, $J_i=1$ for $i=0,1,\ldots, I$.  Their methods consist of finding an approximate critical point of a sequence of penalty subproblems
}
\beq \label{exasubpm0}
\min\limits_{x\in X}  \phi_0(x)+\zeta_0(x) -\psi_0(x) + \rho_k \left[\max\limits_{1\le i \le I} \{\phi_{i}(x)+\zeta_i(x)-\psi_i(x)\}\right]_+,
\eeq
where $\rho_k >0$ is updated by a specific scheme and $[a]_+=\max\{a,0\}$ for any $a\in \Re$.
In particular,{ for DCA1, a DC algorithm is used to find an approximate critical point of \eqref{exasubpm0} by solving a sequence of problems
\beq
x^{l+1} \in \Argmin\limits_{x\in X}  \ u(x) - (s^l)^T x,
\label{convex-subprob2}
\eeq
where $s^l \in \partial v(x^l)$, and
\begin{align}
& u(x) = \phi_0(x)+\zeta_0(x) + \rho_k \max\left\{\max\limits_{1 \le i \le I}\{\phi_{i}(x)+\zeta_i(x)+\sum^I_{j=1,j\neq i} \psi_j(x)\}, \sum^I_{j=1} \psi_j(x)\right\}, \label{fun-p}\\
& v(x) = \psi_0(x) + \rho_k \sum^I_{j=1} \psi_j(x). \label{fun-q}
\end{align}
In addition, for DCA2, a majorization algorithm is used to find an approximate critical point of \eqref{exasubpm0} by solving a sequence of problems
\begin{align}
\begin{array}{rl}
(x^{l+1}, t_{l+1}) \in \Argmin\limits_{x\in X, t\ge 0}& \phi_0(x)+\zeta_0(x)-(s^l_{\psi_0})^T x + \rho_k t \\
\mbox{s.t.} &  \phi_{i}(x)+\zeta_i(x)-\psi_i(x^l)-(s^l_{\psi_i})^T(x-x^l) - t \leq 0, \ \ \forall i=1,\ldots,I, \label{convex-subprob3}
\end{array}
\end{align}
where  $s^l_{\psi_i} \in \partial \psi_i(x^l)$ for $i=0,1,\ldots, I$.
 It was shown in \cite[Theorems 1 and 2]{Thi14a} that any accumulation point of the solution sequence generated by DCA1 and DCA2 is a KKT point that is defined in \eqref{KKT}.  Nevertheless, the proofs of \cite[Theorems 1 and 2]{Thi14a} are based on the assumption that the extended Mangasarian-Fromowitz constraint qualification (EMFCQ) holds at every {\it infeasible accumulation point} and every boundary accumulation point of the solution sequence, which is unreasonable because a constraint qualification is usually assumed to hold at a feasible point rather than infeasible points. Besides this, it is typically difficult to find an {\it exact} solution of subproblems \eqref{convex-subprob2} and
 \eqref{convex-subprob3}, and thus the penalty type of methods \cite{Thi14a} are not practical.}

Recently, Pang et al.\ \cite{Pang} explored the problem structure  and proposed a novel enhanced DCA (EDCA) for solving a special case of (\ref{e1}) with
$I=0$ by solving a number of convex approximation problems per iteration.
They showed that every accumulation point of the solution sequence of EDCA is a
{\it directional-stationary} (D-stationary) {\it point} of the problem.
Besides, Beck and Hallak~\cite{Beck20} proposed a novel feasible descent method for finding a D-stationary point of a class of problems in the form of  $\min\{f(x) - g(x) : x\in X\}$, where $f$ is a continuously differentiable function, $g$ is a convex function, and $X$ is a polyhedral set, which can be applied to a special case of \eqref{e1} with $I=0$, $\zeta_0 \equiv 0$, and $X$ being a polyhedral set.
In addition, assuming that a feasible point of problem \eqref{e1} is available, Pang et al.\ \cite{Pang} proposed an EDCA for solving  (\ref{e1}) by solving a
number of convex approximation problems similar to \eqref{convex-subprob1} per iteration.  They showed that under some suitable constraint qualification, every accumulation point of the generated solution sequence is a {\it Bouligand-stationary} (B-stationary) {\it point} of the problem, which is generally stronger than a usual KKT point.

Although the aforementioned EDCA \cite{Pang} enjoys nice theoretical convergence properties, it is not applicable to problem  \eqref{e1} when a feasible point is not available. To overcome this issue, Pang et al.~\cite{Pang} proposed a penalty approach to solving  \eqref{e1}, which consists of finding an {\it exact} D-stationary point $x^k$ of a sequence of penalty subproblems \eqref{exasubpm0} with $0 < \rho_k \uparrow \infty$. They showed that any feasible accumulation point $x^\infty$ of the sequence $\{x^k\}$ is a B-stationary point of \eqref{e1} if $x^\infty$ satisfies a pointwise Slater constraint qualification (PSCQ) and additionally that the set $\{ j:  \psi_{i,j}(x^\infty) = \psi_i(x^\infty)\}$ is a singleton for $i=1,\ldots, I$.
The latter condition appears to be rather strong because it implies that $\psi_i$ is strictly differentiable at $x^\infty$ for $i=1,\ldots, I$, which generally does not hold, and moreover, the $\psi_i$'s associated with the {\it inactive} constraints are unusually involved. Besides, this penalty approach faces some practical issues. Indeed, as mentioned in \cite{Pang}, problem \eqref{exasubpm0} can be rewritten as
\begin{equation}\label{exasubpm0-1}
\min\limits_{x\in X} \underbrace{u(x)}_{\rm convex} -\underbrace{v(x)}_{\rm convex},
\end{equation}
where $u$ and $v$ are defined in \eqref{fun-p} and \eqref{fun-q}, respectively.
Also, $v$ can be rewritten as
\[
v(x)= \max\left\{ \psi_{0,j_0}(x)+\rho_k\sum^I_{i=1} \psi_{i,j_i}(x) \ \Big| \ 1\le j_i \le J_i,  \quad \forall i=0, 1, \ldots, I \right\}.
\]
It thus follows that \eqref{exasubpm0-1} is a special case of \eqref{e1} with $I=0$. As suggested in \cite{Pang}, problem \eqref{exasubpm0-1} is solved by the aforementioned EDCA, which generates
a sequence any of whose accumulation points is a D-stationary point of \eqref{exasubpm0-1} and hence
of \eqref{exasubpm0}. Therefore, the EDCA is generally only able to produce an {\it approximate} D-stationary point of \eqref{exasubpm0}, but not an {\it exact} one as required by this penalty approach.  In addition, when applied to \eqref{exasubpm0-1}, the EDCA needs to find the exact solution of a number of subproblems in the form of
\begin{equation}\label{exasubpm0-2}
\min\limits_{x\in X} u(x)+\frac{c}{2} \|x-y\|^2
\end{equation}
for some $c>0$ and $y\in\Re^n$, where $u$ is defined in \eqref{fun-p}. Though problem \eqref{exasubpm0-2} is convex, it is typically impossible to find its exact solution due to the sophistication of $u$.
%
%


Motivated by the above points, we propose in this paper a penalty method for solving \eqref{e1} that consists
of a sequence of penalty subproblems in the form of
 \begin{equation}\label{exasubpm}
\min\limits_{x\in X}F_{\rho}(x),
\end{equation}
where
\begin{equation}\label{EXAF}
F_\rho(x)= F(x)+\rho\sum\limits_{i=1}^I\left[\phi_i(x)+\zeta_i(x)-\psi_i(x)\right]_+^p
\end{equation}
with $p\geq1$.\footnote{When $p = 1$, $F_\rho(x)$ has a nonsmooth penalty term and thus may possess a nice exact penalty property.} At each iteration our method only needs an {\it approximate} D-stationary point of the penalty subproblem, which can be efficiently computed by a successive convex approximation method proposed in this paper.
We show that any feasible accumulation point $x^\infty$ of
the solution sequence of our method is a B-stationary point of \eqref{e1} if $x^\infty$ satisfies a PSCQ condition. Compared to the aforementioned convergence result in \cite{Pang}, our result does not require the assumption that the set $\{ j:  \psi_{i,j}(x^\infty) = \psi_i(x^\infty)\}$ is a singleton for $i=1,\ldots, I$. As a consequence, any feasible accumulation point $x^\infty$ of the solution sequence generated by our penalty method can be a B-stationary point of \eqref{e1} even when some of $\psi_i$'s are \emph{non-differentiable} at $x^\infty$. Besides, the PSCQ condition used in our result is generally weaker than that in \cite{Pang}. In fact, we provide an example (see Remark \ref{B-statpt}) for which the PSCQ in our paper holds while the one in \cite{Pang} fails to hold.

In addition, we propose an augmented Lagrangian (AL) method for solving \eqref{e1}, which consists of finding an approximate D-stationary point to a sequence of AL subproblems. Under the same assumptions as those for the penalty method, we show that each accumulation point of the solution sequence is a B-stationary point of \eqref{e1}. We also show that each accumulation point of a set of auxiliary Lagrangian multiplier sequences together with the accumulation point of the solution sequence satisfies a KKT type of optimality condition of \eqref{e1}. Moreover, we provide an example to illustrate the convergence of our AL method.

{We also propose a successive convex approximation method for computing an approximate D-stationary point of the aforementioned penalty and AL subproblems. The proposed method only solves a single convex problem in each iteration, while  the EDCA \cite{Pang} needs to solve a number of convex problems per iteration. It is therefore practically
more  efficient than the latter method.


The rest of the paper is organized as follows. We present some technical preliminaries in Section \ref{tech-prelim} regarding the tangent cone and B-stationary points of problem \eqref{e1}. In Section \ref{penalty-mainprob} we propose a penalty method for solving \eqref{e1} and study its convergence. In Section \ref{alm} we propose an augmented Lagrangian method for solving \eqref{e1} and study its convergence. In Section \ref{penalty-subprob} we propose a successive convex approximation method for solving the penalty and AL subproblems and study  its convergence. We present in Section \ref{Num-Res} some numerical results of the proposed methods.
Finally, in Section \ref{conclude} we make some concluding remarks.}

\subsection{Notation} \label{notation}
Throughout this paper, $\|\cdot\|$ stands for the Euclidean norm and $\Re^n_+$ denotes the nonnegative orthant of the $n$-dimensional real vector space.
We denote $\mathcal{J}_i=\{1, 2, \ldots, J_i\}$ for $i=0,1,\ldots, I$, $\mathcal{I}=\{1,2,\ldots,I\}$, and
$$\mathcal{J}=\{(j_1,j_2,\ldots,j_I)\ |\ j_i\in\mathcal{J}_i,\ \ \forall i\in\mathcal{I}\}.$$
We denote an element of $\mathcal{J}$ by $\mathbb{j}$, i.e., $\mathbb{j}=(j_1,j_2,\ldots,j_I)$
for some $j_i\in\mathcal{J}_i$ for all $i\in\mathcal{I}$.
We use $\Omega$ to denote the feasible region of problem (\ref{e1}), that is,
\[
\Omega= \{x\in X\ |\ \phi_{i}(x)+\zeta_i(x)-\psi_{i}(x)\leq 0,\quad\forall i\in \mathcal{I}\}.
\]
For any $\mathbb{j}=(j_1,j_2,\ldots,j_I)\in \mathcal{J}$, we denote
\begin{equation}\label{feaseq2}
\Omega_{\mathbb{j}}= \{x\in X\ |\ \phi_{i}(x)+\zeta_i(x)-\psi_{i, j_i}(x)\leq 0,\quad  \forall i\in\mathcal{I}\}.
\end{equation}
For any $x\in X$, we denote
\[
\ba{lcl}
\mathcal{J}_i(x) &=& \{j\in\mathcal{J}_i\ |\ \psi_i(x)=\psi_{i,j}(x)\}, \quad \forall i\in\{0,1,\ldots,I\}, \\
\mathcal{J}(x) &=&\{(j_1,j_2,\ldots,j_I)\ |\ j_i\in \mathcal{J}_i(x), \quad \forall i\in \cI\}.
\ea
\]
Clearly, $\mathcal{J}(x) = \mathcal{J}_1(x) \times \cdots \times \mathcal{J}_I(x)$.
For any $x\in X$, let
\[
\ba{l}
\mathcal{I}_>(x)=\{i\in \mathcal{I}\ |\ \phi_i(x)+\zeta_i(x)- \psi_{i}(x)>0\}, \\
\mathcal{I}_=(x)=\{i\in \mathcal{I}\ |\ \phi_i(x)+\zeta_i(x)- \psi_{i}(x)=0\}, \\
\mathcal{I}_<(x)=\{i\in \mathcal{I} \ |\ \phi_i(x)+\zeta_i(x)- \psi_{i}(x)<0\}.
\ea
\]

We now recall some notations from \cite{Rockafellar}. Let $Y\subseteq \Re^n$ and $x \in Y$. The tangent cone of $Y$ at $x$ is denoted by $\mathcal{T}_Y(x)$, i.e.,
$$\mathcal{T}_Y(x)=\left\{d\in\Re^n\Big|\ \exists x^k\in Y,\ x^k\to x,\ \tau_k\downarrow 0 \ {\rm such \ that}\ d=\lim\limits_{k\to \infty}\frac{x^k-x}{\tau_k}\right\}.$$
Also, the normal cone of $Y$ at $x$ is denoted by $\mathcal{N}_Y(x)$. If $Y$ is a closed convex set, $\mathcal{T}_Y(x)$ and $\mathcal{N}_Y(x)$ can be represented as follows:
\begin{align}
\mathcal{T}_Y(x) &= {\rm cl}\left(\{\tau(\bar{x}-x) \ |\ \forall \tau\geq0,\ \forall \bar{x}\in Y\}\right), \label{TC}\\
\mathcal{N}_Y(x) &= \{v\in\Re^n\ |\ v^T(\bar{x}-x)\leq0,\ \forall \bar{x}\in Y\}, \nn
\end{align}
where ${\rm cl}(\cdot)$ is the closure of the associated set.
For a function $f:\Re^n \to \Re\cup \{\infty\}$, the directional derivative of $f$ at a point $x$ in its domain
along a direction $d\in \Re^n$ is defined as
$$f'(x; d)=\lim\limits_{\tau\downarrow 0}\frac{f(x+\tau d)-f(x)}{\tau}.$$
Specifically, by \eqref{e2} and the well-known Danskin's Theorem (e.g., see \cite[Proposition B.25]{Bertsekas}), one can deduce that
\begin{equation}\label{Danskin}
\psi_i'(x;d)=\max\limits_{j\in \mathcal{J}_i(x)} \nabla \psi_{i,j}(x)^Td, \quad \forall x\in X,\ d\in\mathcal{T}_X(x),\ i\in\{0,1,\ldots,I\}.
\end{equation}

A point $x\in Y$ is called a {\em B-stationary point} of $f$ on $Y$ if it satisfies
\begin{equation}\label{B-sp}
f'(x; d)\geq0,\qquad \forall d\in \mathcal{T}_Y(x).
\end{equation}
If $Y$ is a closed convex set and $x\in Y$, $x$ is called a {\em D-stationary point} of $f$ on $Y$ if
\beq \label{D-sp}
f'(x; \bar{x}-x)\geq0,\qquad \forall \bar{x} \in Y.
\eeq
It follows from \eqref{TC}, \eqref{B-sp} and \eqref{D-sp} that a B-stationary point of $f$ on $Y$ reduces to
a D-stationary point  when $Y$ is a closed convex set. See \cite{Pang07,Pang} for more discussion.

For a smooth function $f$ on $X$ and $\bar{x}\in X$, we define
\beq \label{lg}
\ell_f(x; \bx) = f(\bx) + \nabla f(\bx)^T(x-\bx),
\eeq
which is the linearization of $f$ at $\bar x$. Clearly, $f(x)\geq\ell_f(x; \bx)$ when $f$ is convex on $X$.

\section{Technical preliminaries} \label{tech-prelim}

Due to the nonsmoothness and sophistication of the constraints of \eqref{e1}, it is
generally difficult to characterize the tangent cone $\mathcal{T}_{\Omega}(\bar x)$ at a
point $\bar x\in\Omega$, where $\Omega$ is the feasible region of \eqref{e1}. In this
section,  we provide some characterization of  $\mathcal{T}_{\Omega}(\bar x)$ under the
PSCQ condition at $\bar x$ by exploiting the special structure of $\psi_i(\bar x)$, which
is a generalization of a result by Pang et al.\ \cite{Pang} for a special case with $I=1$
and $\phi_1 \equiv 0$.  As a consequence, we provide a characterization for a B-stationary point of \eqref{e1}. In addition, under some suitable assumption we also provide a KKT type
of characterization for a B-stationary point of \eqref{e1}.

{One can easily observe from \eqref{e1}, \eqref{e2} and \eqref{feaseq2} that  $\Omega=\bigcup_{\mathbb{j}\in \mathcal{J}}\Omega_{\mathbb{j}}$.} In addition, it is not hard to observe that
\[
\mathcal{T}_{\Omega}(\bar{x})=\mathcal{T}_{\bigcup_{\mathbb{j}\in \mathcal{J}(\bar{x})}\Omega_{\mathbb{j}}}(\bar{x}), \quad \forall \bar{x}\in \Omega.
\]
It follows from this relation and \cite[Table 4.1]{Aubin} that
\beq \label{propfea}
\mathcal{T}_{\Omega}(\bar{x})=\bigcup_{\mathbb{j}\in \mathcal{J}(\bar{x})}\mathcal{T}_{\Omega_{\mathbb{j}}}(\bar{x}), \quad \forall \bar{x}\in \Omega.
\eeq

From \eqref{propfea}, one can see that to characterize $\mathcal{T}_{\Omega}(\bar{x})$ at a point $\bar x\in\Omega$, it suffices to characterize $\mathcal{T}_{\Omega_{\mathbb{j}}}(\bar{x})$ for every $\mathbb{j}\in\mathcal{J}(\bar{x})$. To proceed, let $\bar{x}\in \Omega$ and $\mathbb{j}=(j_1,\ldots,j_I)\in \mathcal{J}(\bar{x})$. In a similar vein as in \cite{Pang}, we define
\begin{align}
Y_{\mathbb{j}}(\bar{x}) & =\{x\in X\ |\ \phi_i(x)+\zeta_i(x)-\ell_{\psi_{i,j_i}}(x;\bar{x}) \leq 0,\quad\forall i\in\mathcal{I}\}, \label{cvl} \\
C_{\mathbb{j}}(\bar{x}) & = \{d\in \mathcal{T}_{X}(\bar{x})\ |\ \nabla \phi_i(\bar{x})^Td+\zeta'_i(\bar{x};d)-\nabla \psi_{i,j_i}(\bar{x})^Td\leq0,\quad \forall i\in \mathcal{I}_=(\bar{x})\}. \label{lin-cone}
\end{align}
We are now ready to provide a characterization of $\mathcal{T}_{\Omega_{\mathbb{j}}}(\bar{x})$, whose proof is similar to that of Propositions 2 and 3 in \cite{Pang} and thus omitted.

\begin{proposition}\label{consqua}
Let $\bar{x}\in\Omega$ be such that $\cI_=(\bx)\neq\emptyset$ and $\mathbb{j}=(j_1,\ldots, j_I)\in \mathcal{J}(\bar{x})$. Suppose that there exists some $\bar{d}\in\mathcal{T}_X(\bar{x})$ such that
\begin{equation}\label{px}
\nabla \phi_i(\bar{x})^T\bar{d}+\zeta_i'(\bar{x};\bar{d})-\nabla \psi_{i,j_i}(\bar{x})^T\bar{d}<0, \qquad \forall i \in \mathcal{I}_=(\bar{x}),
\end{equation}
or that $X$ is a polyhedral set, $\phi_{i}$ is affine and $\zeta_i$ is piecewise affine on $X$ for every $i\in\mathcal{I}_=(\bar{x})$. Then
\begin{equation*}\label{eq3}
\mathcal{T}_{Y_{\mathbb{j}}(\bar{x})}(\bar{x})=\mathcal{T}_{\Omega_{\mathbb{j}}}(\bar{x})={C_{\mathbb{j}}(\bar{x})}.
\end{equation*}
\end{proposition}

From Proposition \ref{consqua}, we see that condition \eqref{px} is crucial for characterizing $\mathcal{T}_{\Omega_{\mathbb{j}}}(\bar{x})$ for a given $\mathbb{j}\in \mathcal{J}(\bar{x})$. In view of this and $\mathcal{T}_{\Omega}(\bar{x})=\bigcup_{\mathbb{j}\in \mathcal{J}(\bar{x})}\mathcal{T}_{\Omega_{\mathbb{j}}}(\bar{x})$, it is natural to introduce the following condition under which a characterization of $\mathcal{T}_{\Omega}(\bar{x})$ can be obtained.

\begin{definition}\label{slater}
Let $\bar{x}\in\Omega$ be such that $\mathcal{I}_=(\bar{x})\neq\emptyset$. The pointwise Slater constraint qualification {\rm (PSCQ)} is said to hold for the set $\Omega$ at $\bar{x}$ if for every $\mathbb{j}=(j_1,\ldots, j_I)\in \mathcal{J}(\bar{x})$, there exists some $d_\mathbb{j}\in\mathcal{T}_X(\bar{x})$ such that
\begin{equation}
\label{eq:slater}
\nabla \phi_i(\bar{x})^Td_\mathbb{j}+\zeta_i'(\bar{x};d_\mathbb{j})-\nabla \psi_{i,j_i}(\bar{x})^Td_\mathbb{j}<0, \qquad\forall i\in\mathcal{I}_=(\bar{x}).
\end{equation}
\end{definition}

\begin{remark} \label{B-statpt}
(i) It shall be mentioned that the concept of PSCQ is always associated with a specific algebraic representation of the underlying set. Strictly speaking,  the PSCQ for the set $\Omega$ in Definition \ref{slater} is based on its algebraic representation given by
\[
\Omega = \left\{x\in X \ \left| \ \phi_{i}(x)+\zeta_i(x)-\max\limits_{1\leq j\leq J_i}\{\psi_{i,j}(x)\}\leq 0, \quad\forall i=1,\ldots,I.\right.\right\}.
\]
Note that the set $\Omega$ also admits the following equivalent algebraic representation:
\begin{equation}\label{new-repre-omega}
\Omega = \left\{x\in X \ \left|  \ \max_{1\le i \le I} \left\{\phi_{i}(x)+\zeta_i(x)+\sum^I_{j=1,j\neq i} \psi_j(x) \right\}- \max\limits_{(j_1,\ldots, j_I)\in\cJ} \sum^I_{i=1}\psi_{i,j_i}(x) \leq 0 \right.\right\}.
\end{equation}
Such a representation of $\Omega$ is used in \cite{Pang} for reformulating problem \eqref{e1} with more than one DC constraints into an equivalent problem with a single DC constraint. While it provides a simplified treatment from a theoretical point of view, the PSCQ defined in \cite{Pang} by using \eqref{new-repre-omega} is generally stronger than that in Definition \ref{slater}.
In fact, the PSCQ defined in \cite{Pang} by using \eqref{new-repre-omega} says that PSCQ holds at
$\bar{x}\in\Omega$ if for every $\mathbb{j}=(j_1,\ldots, j_I)\in \mathcal{J}(\bar{x})$, there exists some $d_\mathbb{j}\in\mathcal{T}_X(\bar{x})$ such that
\[
\nabla \phi_i(\bx)^Td_\mathbb{j} +\zeta_i'(\bx;d_\mathbb{j})-\nabla \psi_{i,j_i}(\bx)^T d_\mathbb{j}+\sum^I_{\ell=1,\ell\neq i} [\psi'_\ell(\bx;d_\mathbb{j})-\nabla \psi_{\ell,j_\ell}(\bx)^Td_\mathbb{j}] <0,  \quad \forall i\in\mathcal{I}_=(\bx).
\]
One can observe that such PSCQ is generally stronger than our PSCQ in Definition \ref{slater}. That is, if  such PSCQ holds at $\bar{x}\in\Omega$, our PSCQ must also hold at $\bar{x}$, while the converse may not hold. As a counterexample, consider the set
\[
\Omega = \{ x\in\Re^2\ | \ -(x_1^2 + x_1 + x_2)\leq 0, \ -\max\{x_1^2 - x_1 + x_2, \ 2x_2\} \leq 0 \}.
\]
Clearly, it is a special case of the feasible region of \eqref{e1} with $I=2$, $X = \Re^2, \phi_1 \equiv \zeta_1 \equiv \phi_2 \equiv \zeta_2 \equiv 0$, $\psi_1(x) = \psi_{1,1}(x)$, and $\psi_2(x) = \max\{\psi_{2,1}(x), \psi_{2,2}(x)\}$, where $\psi_{1,1}(x) = x_1^2 + x_1 + x_2$, $\psi_{2,1}(x) = x_1^2 - x_1 + x_2$, and $\psi_{2,2}(x) = 2x_2$. Let $\bar{x} = (0,0)^T$. One can verify that for such $\Omega$, our PSCQ holds at $\bar x$, but the PSCQ defined in \cite{Pang} by using \eqref{new-repre-omega} fails to hold at $\bar x$.

(ii) It can be shown that the PSCQ holds for the set $\Omega$ at a feasible point $\bar{x}$ if and only if for every $\mathbb{j}=(j_1,\ldots, j_I)\in \mathcal{J}(\bar{x})$, there exists a Slater point in the set
\begin{equation*}
\left\{x\in X| \phi_i(x)+\zeta_i(x)- \psi_{i,j_i}(\bar{x})-\nabla \psi_{i,j_i}(\bar{x})^T(x-\bar{x}) \le 0, \ \forall i\in\mathcal{I}_=(\bar{x})\right\},
\end{equation*}
which is sometimes more checkable than the conditions in \eqref{eq:slater}. From this, one can see that PSCQ is indeed a generalization of the classical Slater's condition.
\end{remark}

As a consequence of \eqref{propfea}, Definition \ref{slater} and Proposition \ref{consqua}, we can obtain the following characterization of $\mathcal{T}_{\Omega}(\bar{x})$ at a point $\bar{x}\in\Omega$ with $\cI_=(\bx)\neq\emptyset$.

\begin{corollary}\label{cor1}
Let $\bar{x}\in\Omega$ be such that $\cI_=(\bx)\neq\emptyset$. Suppose that the PSCQ holds for $\Omega$ at $\bar{x}$ or that $X$ is a polyhedral set, $\phi_{i}$ is affine and $\zeta_i$ is piecewise affine on $X$ for every $i\in\mathcal{I}_=(\bar{x})$. Then
\begin{equation*}\label{equility}
\mathcal{T}_{\Omega}(\bar{x})=\bigcup_{\mathbb{j}\in\mathcal{J}(\bar{x})}C_\mathbb{j}(\bar{x})=\bigcup_{\mathbb{j}\in\mathcal{J}(\bar{x})}\mathcal{T}_{Y_{\mathbb{j}}(\bar{x})}(\bar{x}).
\end{equation*}
\end{corollary}

From Corollary \ref{cor1}, we immediately obtain the following characterization of a B-stationary point of \eqref{e1}.

\begin{theorem}
\label{suff-cond-B}
Let $\bar{x}\in\Omega$ be such that $\mathcal{I}_=(\bar{x})\neq\emptyset$. Suppose that the  {\rm PSCQ} holds for $\Omega$ at $\bar{x}$, or that $X$ is a polyhedral set, $\phi_{i}$ is affine and $\zeta_i$ is piecewise affine on $X$ for every $i\in\mathcal{I}_=(\bar{x})$.  Then $\bx$ is a B-stationary point of problem \eqref{e1} if and only if $F'(\bx; d) \ge 0$ for all
$d \in \mathcal{T}_{Y_{\mathbb{j}}(\bar{x})}(\bar{x})$ and $\mathbb{j}\in\mathcal{J}(\bar{x})$.
\end{theorem}


Before ending this section, we provide a KKT type of characterization of a B-stationary point of problem (\ref{e1}), whose proof is given in Appendix \ref{proof-KKT}.

\begin{theorem}\label{thm:KKT}
Let $\bx\in\Omega$ be such that $\cI_=(\bx)\neq\emptyset$. Suppose that the  {\rm PSCQ} holds for $\Omega$ at $\bar{x}$, or that $X$ is a polyhedral set, $\phi_{i}$ is affine and $\zeta_i$ is piecewise affine on $X$ for every $i\in\mathcal{I}_=(\bar{x})$. Then $\bx$ is a B-stationary point of problem \eqref{e1} if and only if for every $j_0\in\cJ_0(\bx)$ and every $\bj=(j_1,\cdots,j_I)\in\cJ(\bx)$, there exists a vector of Lagrangian multipliers $\lambda^{j_0,\bj}=(\lambda_1^{j_0,\bj},\ldots,\lambda^{j_0,\bj}_I)$ satisfying that
\beqa
&& \lambda_i^{j_0,\bj}\geq0,
\quad \lambda_i^{j_0,\bj}\left[\phi_i(\bx)+\zeta_i(\bx)-\psi_{i,j_i}(\bx)\right]=0,\quad\forall i\in\cI, \label{kkt-1} \\
&& 0\in \nabla\phi_0(\bx)+\partial \zeta_0(\bx)-\nabla\psi_{0,j_0}(\bx)+\sum\limits_{i=1}^I\lambda_i^{j_0,\bj}\left[\nabla\phi_i(\bx)+\partial\zeta_i(\bx)-\nabla\psi_{i,j_i}(\bx)\right]+\mathcal{N}_X(\bx).\qquad \label{kkt-2}
\eeqa
\end{theorem}

\section{A penalty method for DC program \eqref{e1}}
\label{penalty-mainprob}
In this section we propose a penalty method for solving problem \eqref{e1}, which consists of finding an approximate solution to a sequence of penalty subproblems in the form of \eqref{exasubpm}. Before proceeding, we introduce some notations that will be used shortly.

Let $\epsilon>0$ be given. Define
\beqa
&&\mathcal{J}_{0,\epsilon}(x)=\{j\in\mathcal{J}_0 \ |\ \psi_0(x)\le \psi_{0,j}(x)+\epsilon\}, \label{J0-index}\\
&&\mathcal{J}_{\epsilon}(x)=\{(j_1,\ldots, j_I)\in\mathcal{J} \ |\ \psi_i(x) \le \psi_{i,j_i}(x)+\epsilon,\ \forall i\in\mathcal{I}\}.\label{J-index}
\eeqa
Given any $\bar{x}\in X$, let
\beq \label{hphi}
\hphi_i(x;\bar{x}) = \phi_i(\bar{x}) + \nabla \phi_i(\bar{x})^T (x-\bar x) + L_i \|x-\bar x\|^2/2, \quad \forall i\in\{0,1,\ldots, I\},
\eeq
where $L_i$ is the Lipschitz constant associated with $\nabla\phi_i$ on $X$. Moreover, for every $i\in\cI$, we let $L_i=0$ if $\phi_i$ is affine.
For any $\bar{x}\in X$, $p\ge 1$, $j_0\in \mathcal{J}_0$ and $\mathbb{j}=(j_1,\ldots,j_I) \in \mathcal{J}$, we define
\begin{eqnarray}
&&F_\rho(x,j_0,\mathbb{j})= \phi_0(x)+\zeta_0(x)-\psi_{0,j_0}(x)+\rho\sum\limits_{i=1}^I[\phi_i(x)+\zeta_i(x)-\psi_{i,j_i}(x)]_+^p, \label{exasmsbp}\\
&& Q_\rho(x;\bar{x},j_0,\mathbb{j}) = \hphi_0(x;\bar{x}) + \zeta_0(x) - \ell_{\psi_{0,j_0}}(x; \bx) + \rho\sum\limits_{i=1}^I[\hphi_i(x;\bar{x}) + \zeta_i(x) - \ell_{\psi_{i,j_i}}(x; \bx)]_+^p. \qquad
\label{exaqappfun}
\end{eqnarray}

\begin{remark}
(i) For any $\bx \in X$, $\mathcal{J}_{0}(\bx) \subseteq \mathcal{J}_{0,\epsilon}(x)$
and $\mathcal{J}(\bx) \subseteq \mathcal{J}_{\epsilon}(x)$ for all $x \in X$ sufficiently close to $\bx$.

(ii) By the Lipschitz continuity of $\nabla \phi_i$, one has $\phi_i(x) \le \hphi_i(x;\bar{x})$ for all $x\in X$. Also, since $\psi_{i,j_i}$ is convex on $X$, $\psi_{i,j_i}(x) \ge \ell_{\psi_{i,j_i}}(x; \bx)$ for all $x\in X$ and $\bx\in X$.
In view of  these, \eqref{exasmsbp} and \eqref{exaqappfun}, we can observe that $F_\rho(x,j_0,\mathbb{j})
\le Q_\rho(x;\bar{x},j_0,\mathbb{j})$ for all $x\in X$.
\end{remark}


We now present a penalty method for solving problem \eqref{e1} and establish its convergence. The details of the penalty method are presented as follows.

\begin{algorithm}
\label{PENM}
\normalfont
\mbox{}
\begin{itemize}
\item[0.] Input $\epsilon>0$, $\rho_0>0$, $\sigma>1$, and a sequence $\{\eta_k\}\subset\Re_+$ such that $\eta_k\to0$. Set $k\leftarrow 0$.
\item[1.] Find an approximate solution $x^k$ of the penalty subproblem
\begin{equation}\label{sbpp}
\min\limits_{x\in X}F_{\rho_k}(x)
\end{equation}
such that $x^k \in X$ and
\begin{equation}\label{exawww10}
F_{\rho_k}(x^{k})\leq Q_{\rho_k}(x;x^{k},j_0,\mathbb{j})+\eta_k,\qquad\forall x\in X
\end{equation}
for every $j_0\in\mathcal{J}_{0,\epsilon}(x^{k})$ and $\mathbb{j}\in\mathcal{J}_{\epsilon}(x^{k})$, where
$F_{\rho_k}$ is defined in \eqref{EXAF}.
\item[2.] Set $\rho_{k+1} \leftarrow \sigma \rho_k$.
\item[3.] Set $k \leftarrow k+1$, and go to Step 1.
\end{itemize}
{\bf End.}
\end{algorithm}


To make the above penalty method complete, we need to address how to find an approximate solution $x^k\in X$ for subproblem (\ref{sbpp}) satisfying (\ref{exawww10}) as required in Step 1. We will leave this discussion in Section \ref{penalty-subprob}. For the time being, we establish the main convergence result regarding this method for solving problem (\ref{e1}).

\begin{theorem}\label{inexctct}
 Let $\{x^k\}$ be generated by Algorithm \ref{PENM}. Assume that $\{x^k\}_{k\in \mathcal{K}}$ converges to $x^\infty$ for some subsequence $\cK$. Then the following statements hold.
\bi
\item[(i)] $x^\infty$ is a D-stationary point of the problem
\begin{equation}\label{exaaaa333}
\min\limits_{x\in X}\sum\limits_{i=1}^I\underbrace{[\phi_i(x)+\zeta_i(x)-\psi_{i}(x)]_+^p}_{h_i(x)}.
\end{equation}
\item[(ii)] If $x^{\infty}\in\Omega$ and $\mathcal{I}_=(x^\infty)=\emptyset$, then $x^\infty$ is a D-stationary point of the problem
\begin{equation}\label{ee3}
\min\limits_{x\in X} F(x).
\end{equation}
\item[(iii)] If $x^{\infty}\in\Omega$, $\mathcal{I}_=(x^\infty)\neq\emptyset$, and moreover, the PSCQ holds for $\Omega$ at $x^\infty$, then $x^\infty$ is a B-stationary point of problem (\ref{e1}).
\item[(iv)] If $x^{\infty}\in\Omega$, $\mathcal{I}_=(x^\infty)\neq\emptyset$, $X$ is a polyhedral set, and moreover, for every $i\in\mathcal{I}_=(x^{\infty})$ and $j_i\in\cJ_i(x^\infty)$, $\phi_{i}$ and $\psi_{i,j_i}$ are affine and $\zeta_i$ is piecewise affine on $X$, then $x^\infty$ is a B-stationary point of problem (\ref{e1}).
\ei
\end{theorem}

\begin{proof}  Since $\{x^k\}_{k\in \cK}$ converges to $x^\infty$,  one has $\mathcal{J}_0(x^\infty)\subseteq \mathcal{J}_{0,\epsilon}(x^k)$ and $\mathcal{J}(x^\infty)\subseteq \mathcal{J}_\epsilon(x^k)$ for sufficiently large $k\in \cK$. It thus follows from (\ref{EXAF}), (\ref{exaqappfun}) and (\ref{exawww10}) that for $k\in \cK$ sufficiently large, we have
\begin{equation}\label{exaaaa111}
\begin{array}{ll}
\displaystyle\phi_0(x^k)+\zeta_0(x^k)-\psi_0(x^k)+\rho_k\sum\limits_{i=1}^I[\phi_i(x^k)+\zeta_i(x^k)-\psi_i(x^k)]_+^p-\eta_k\\
\displaystyle \leq \hat{\phi}_0(x;x^k)+\zeta_0(x)-\ell_{\psi_{0,j_0}}(x;x^k)+\rho_k\sum\limits_{i=1}^I[\hat{\phi}_i(x;x^k)+\zeta_i(x)-\ell_{\psi_{i,j_i}}(x;x^k)]_+^{p}
\end{array}
\end{equation}
for all $j_0\in\mathcal{J}_0(x^\infty)$, $\mathbb{j}=(j_1,\ldots,j_I)\in\mathcal{J}(x^\infty)$ and $x\in X$.

(i) {{In order to prove statement (i), we first show that
\beq \label{eq-thm3.1(i)}
x^\infty \in \Argmin_{x\in X}\sum\limits_{i=1}^I\underbrace{\left[\hat{\phi}_i(x;x^\infty)+\zeta_i(x)-\ell_{\psi_{i,j_i}}(x;x^\infty)\right]_+^p}_{\hh_{i,j_i}(x)}, \quad \forall \bj=(j_1,\ldots, j_I) \in \mathcal{J}(x^\infty).
\eeq
Indeed, notice that}} $\{x^k\}_{k\in \cK}\to x^\infty$ and $\{\rho_k\} \to \infty$. Dividing both sides of (\ref{exaaaa111}) by $\rho_k$  and
taking limits as $\cK\ni k \to \infty$ yield
\beq \label{constr-opt}
\sum\limits_{i=1}^I\left[\phi_i(x^\infty)+\zeta_i(x^\infty)-\psi_i(x^\infty)\right]_+^p \le \sum\limits_{i=1}^I\left[\hat{\phi}_i(x;x^\infty)+\zeta_i(x)-\ell_{\psi_{i,j_i}}(x;x^\infty)\right]_+^p
\eeq
for any $x\in X$ and $\mathbb{j}=(j_1,\ldots,j_I) \in \mathcal{J}(x^\infty)$. In view of \eqref{lg} and \eqref{hphi}, one can observe that
\beq \label{value-xinfty}
\hat{\phi}_i(x^\infty;x^\infty)=\phi_i(x^\infty), \quad \ell_{\psi_{i,j_i}}(x^\infty;x^\infty) = \psi_{i,j_i}(x^\infty)=
\psi_i(x^\infty)
\eeq
for any $i \in {\cal I}$ and $\mathbb{j}=(j_1,\ldots,j_I) \in \mathcal{J}(x^\infty)$. It follows from \eqref{constr-opt} and \eqref{value-xinfty} that {{\eqref{eq-thm3.1(i)} holds.

We are now ready to complete the proof of statement (i). Clearly,  the relation \eqref{eq-thm3.1(i)}}} yields $\hh'_{\bj}(x^\infty; d) \ge 0$ for all $d\in \mathcal{T}_X(x^\infty)$, where $\hh_\bj(x)=\sum\limits_{i=1}^I \hh_{i,j_i}(x)$ for every $\bj=(j_1,\ldots, j_I) \in \mathcal{J}(x^\infty)$. By virtue of
\eqref{lg} and \eqref{hphi}, it is not hard to verify that for every $\bj=(j_1,\ldots, j_I) \in \mathcal{J}(x^\infty)$,
\[
\hh'_{i,j_i}(x^\infty; d) = \left\{\ba{ll}
\theta_i(x^\infty)[\nabla\phi_i(x^\infty)^T d +\zeta'_i(x^\infty; d)-\nabla \psi_{i,j_i}(x^\infty)^Td], & \mbox{if} \ i \in \cI_>(x^\infty), \\
\theta_i(x^\infty)\left[\nabla\phi_i(x^\infty)^T d +\zeta'_i(x^\infty; d)-\nabla \psi_{i,j_i}(x^\infty)^Td)\right]_+, & \mbox{if} \ i \in \cI_=(x^\infty), \\
0, & \mbox{if} \ i \in \cI_<(x^\infty),
\ea\right.
\]
where $\theta_i(x^\infty)=p[\phi_i(x^\infty)+\zeta_i(x^\infty)-\psi_i(x^\infty)]_+^{p-1}$ for which we assume $0^0=1$.
 In addition, let $h_i$ be defined as in \eqref{exaaaa333} and $h(x)=\sum^I_{i=1}h_i(x)$. By \eqref{Danskin} and \eqref{exaaaa333}, one can observe that
\begin{equation*}
h'_i(x^\infty; d) =\left\{\ba{ll}
\theta_i(x^\infty)[\nabla\phi_i(x^\infty)^T d +\zeta'_i(x^\infty; d) - \max\limits_{j\in \cJ_i(x^\infty)} \nabla \psi_{i,j}(x^\infty)^Td], & \mbox{if} \ i \in \cI_>(x^\infty), \\
\theta_i(x^\infty)[\nabla\phi_i(x^\infty)^T d + \zeta'_i(x^\infty; d)-\max\limits_{j\in \cJ_i(x^\infty)} \nabla \psi_{i,j}(x^\infty)^Td)]_+, & \mbox{if} \ i \in \cI_=(x^\infty), \\
0, & \mbox{if} \ i \in \cI_<(x^\infty).
\ea\right.
\end{equation*}
Hence, for every $d\in \mathcal{T}_X(x^\infty)$, there exists some $\mathbb{j}=(j_1,\ldots,j_I)\in\cJ(x^\infty)$ such that $h'(x^\infty;d)=\hh_\bj'(x^\infty;d)$, which along with $\hh'_{\bj}(x^\infty; d) \ge 0$ implies $h'(x^\infty; d) \ge 0$. It follows from this, $h(x)=\sum^I_{i=1}h_i(x)$ and \eqref{exaaaa333} that statement (i) holds.

(ii) Suppose that $x^{\infty}\in\Omega$ and $\mathcal{I}_=(x^\infty)=\emptyset$. Then
$\mathcal{I}_<(x^\infty)=\mathcal{I}$. {{In order to prove statement (ii), we first show that
\beq \label{eq-thm3.1(ii)}
x^\infty \in \Argmin_{x\in X} \underbrace{\hphi_0(x;x^\infty) + \zeta_0(x)-\ell_{\psi_{0,j_0}}(x;x^\infty)}_{\hh_{0,j_0}(x)},
\qquad \forall  j_0\in\mathcal{J}_0(x^\infty).
\eeq
Indeed, it follows from $\mathcal{I}_<(x^\infty)=\mathcal{I}$ and }}
 \eqref{value-xinfty}  that
\beq \label{strict-feas}
\hat{\phi}_i(x^\infty;x^\infty)+\zeta_i(x^\infty)-\ell_{\psi_{i,j_i}}(x^\infty;x^\infty)=\phi_i(x^\infty)+\zeta_i(x^\infty)-\psi_i(x^\infty)<0
\eeq
for any $i\in\mathcal{I}$ and $j_i \in \cJ_i(x^\infty)$. Let $x\in X$ be arbitrarily chosen. By \eqref{strict-feas}, the continuity of $\zeta_i$ on $X$, and the continuity of $\hphi_i(\cdot;\cdot)$ and $\ell_{\psi_{i,j_i}}(\cdot;\cdot)$ on $X \times X$, one has
\begin{equation}\label{ee1}
\hat{\phi}_i(x^\infty+ t (x-x^\infty);x^k)+\zeta_i(x^\infty+t (x-x^\infty))-\ell_{\psi_{i,j_i}}(x^\infty+t (x-x^\infty);x^k)<0
\end{equation}
for any $j_i\in \mathcal{J}_i(x^\infty)$, $k\in \cK$ sufficiently large and $t>0$ sufficiently small.
Replacing $x$ by $x^\infty+t (x-x^\infty)$ in (\ref{exaaaa111}) and using \eqref{ee1}, we obtain that for every $j_0\in\mathcal{J}_0(x^\infty)$,
\[\phi_0(x^k)+\zeta_0(x^k)-\psi_0(x^k)-\eta_k\leq \hat{\phi}_0(x^\infty+t (x-x^\infty);x^k)+\zeta_0(x^\infty+t (x-x^\infty))-\ell_{\psi_{0,j_0}}(x^\infty+t (x-x^\infty);x^k)
\]
holds for $k\in \cK$ sufficiently large and $t>0$ sufficiently small. Taking limit on both sides of this inequality
as $\cK\ni k \to \infty$, one has that  for every $j_0\in\mathcal{J}_0(x^\infty)$, it holds
\beq \label{lim-ineq}
\phi_0(x^\infty)+\zeta_0(x^\infty)-\psi_0(x^\infty)
\leq \hat{\phi}_0(x^\infty+t (x-x^\infty);x^\infty)+\zeta_0(x^\infty+t (x-x^\infty))-\ell_{\psi_{0,j_0}}(x^\infty+t (x-x^\infty);x^\infty)
\eeq
for $t>0$ sufficiently small. By \eqref{lg}, the linearity of $\ell_{\psi_{0,j_0}}(\cdot;x^\infty)$, the convexity of $\zeta_0$, and $\hphi_0(\cdot; x^\infty)$,
$\hphi_0(x^\infty; x^\infty)=\phi_0(x^\infty)$ and $\ell_{\psi_{0,j_0}}(x^\infty;x^\infty)=\psi_0(x^\infty)$ for $j_0\in\mathcal{J}_0(x^\infty)$, we have that  for $t>0$ sufficiently small,
\[
\ba{ll}
\hat{\phi}_0(x^\infty+t (x-x^\infty);x^\infty)+\zeta_0(x^\infty+t (x-x^\infty))-\ell_{\psi_{0,j_0}}(x^\infty+t (x-x^\infty);x^\infty) \\ [5pt]
\displaystyle \le t [\hphi_0(x;x^\infty) + \zeta_0(x)-\ell_{\psi_{0,j_0}}(x;x^\infty)]+ (1-t) [\phi_0(x^\infty) + \zeta_0(x^\infty)-\psi_0(x^\infty)].
\ea
\]
It follows from this and \eqref{lim-ineq} that for $t>0$ sufficiently small,
\[\phi_0(x^\infty)+\zeta_0(x^\infty)-\psi_0(x^\infty)\le t [\hphi_0(x;x^\infty) + \zeta_0(x)-\ell_{\psi_{0,j_0}}(x;x^\infty)]+ (1-t) [\phi_0(x^\infty) + \zeta_0(x^\infty)-\psi_0(x^\infty)],
\]
which implies that
\[
\phi_0(x^\infty)+\zeta_0(x^\infty)-\psi_0(x^\infty) \le \hphi_0(x;x^\infty) + \zeta_0(x)-\ell_{\psi_{0,j_0}}(x;x^\infty),
\quad \forall x \in X, \ j_0\in\mathcal{J}_0(x^\infty).
\]
By this, $\hphi_0(x^\infty; x^\infty)=\phi_0(x^\infty)$ and $\ell_{\psi_{0,j_0}}(x^\infty;x^\infty)=\psi_0(x^\infty)$, one can see that {{\eqref{eq-thm3.1(ii)} holds.

We are now ready to complete the proof of statement (ii). Indeed, it follows from \eqref{eq-thm3.1(ii)}}} that $\hh_{0,j_0}'(x^\infty;d) \ge 0$ for all $d\in \mathcal{T}_X(x^\infty)$ and $j_0\in\mathcal{J}_0(x^\infty)$, which along with \eqref{lg} and \eqref{hphi} implies that
\[
\nabla \phi_0(x^\infty)^Td + \zeta'_0(x^\infty;d) - \nabla \psi_{0,j_0}(x^\infty)^Td \ge 0, \qquad \forall d\in \mathcal{T}_X(x^\infty), \ j_0\in\mathcal{J}_0(x^\infty).
\]
By this and \eqref{Danskin}, one has $F'(x^\infty;d) \ge 0$ for all $d\in \mathcal{T}_X(x^\infty)$. Hence,
 $x^\infty$ is a D-stationary point of problem (\ref{ee3}).

(iii) Suppose that $x^{\infty}\in\Omega$, $\mathcal{I}_=(x^\infty)\neq\emptyset$, and moreover,  the PSCQ holds for $\Omega$ at $x^\infty$. {{In order to prove statement (iii), we first show that for any $\mathbb{j}=(j_1,\ldots,j_I)\in\mathcal{J}(x^\infty)$, there exists some $\hat x\in Y_\mathbb{j}(x^\infty)$ such that
\beq \label{slater-cond3}
\phi_i(\hat x)+\zeta_i(\hat x)-\ell_{\psi_{i,j_i}}(\hat x; x^\infty)< 0,\quad\forall i\in\mathcal{I}_=(x^\infty).
\eeq
Indeed, let $\mathbb{j}=(j_1,\ldots,j_I)\in\mathcal{J}(x^\infty)$ be arbitrarily chosen. Since PSCQ holds for $\Omega$ at $x^\infty$,}} there exists some $d\in\mathcal{T}_X(x^\infty)$ such that
\beq \label{slater-cond1}
\nabla \phi_i(x^\infty)^Td+\zeta_i'(x^\infty;d)-\nabla \psi_{i,j_i}(x^\infty)^Td<0,\qquad\forall i\in\mathcal{I}_=(x^\infty).
\eeq
Hence, there exist $\{{\hat x}^l\} \subseteq X$ and $\{\alpha_l\} \downarrow 0$ such that $d=\lim_{l\to\infty} ({\hat x}^l-x^\infty)/\alpha_l$, which implies
${\hat x}^l=x^\infty + \alpha_l d +o(\alpha_l)$. It follows from this and \eqref{slater-cond1} that
\[
\ba{ll}
\displaystyle\lim\limits_{l\to\infty} \frac{1}{\alpha_l} \left\{[\phi_i({\hat x}^l)+\zeta_i({\hat x}^l)-\ell_{\psi_{i,j_i}}({\hat x}^l;x^\infty)] -[\phi_i(x^\infty)+\zeta_i(x^\infty)-\psi_{i,j_i}(x^\infty)]\right\} \\ [8pt]
\displaystyle =\nabla \phi_i(x^\infty)^Td+\zeta_i'(x^\infty;d)-\nabla \psi_{i,j_i}(x^\infty)^Td < 0, \qquad\forall i\in\mathcal{I}_=(x^\infty).
\ea
\]
By this and $\phi_i(x^\infty)+\zeta_i(x^\infty)-\psi_{i,j_i}(x^\infty)=0$ for $i\in\mathcal{I}_=(x^\infty)$,  we have that
for $l$ sufficiently large,
\[ \phi_i({\hat x}^l)+\zeta_i({\hat x}^l)-\ell_{\psi_{i,j_i}}({\hat x}^l;x^\infty)<0, \quad\forall i\in\mathcal{I}_=(x^\infty).\]
Therefore, there exists some $\hat x\in Y_\mathbb{j}(x^\infty)$ such that \eqref{slater-cond3} holds.

{We next show that for any $\tau\in(0,1]$, $\mathbb{j}=(j_1,\ldots,j_I)\in\mathcal{J}(x^\infty)$, and $x\in Y_\mathbb{j}(x^\infty)$, there exists some $\hat{t}\in(0,1)$ such that
$x(\hat{t},\tau)\in Y_\mathbb{j}(x^\infty)$ and
\beq  \label{uppbd-xtau3}
\hphi_i(x(\hat t,\tau);x^\infty)+\zeta_i(x(\hat t,\tau)) -\ell_{\psi_{i,j_i}}(x(\hat t,\tau);x^\infty) <0, \quad \forall i \in \cI,
\eeq
where
\beq\label{exahhh333}
x(t,\tau)= x^\infty+t(x(\tau)-x^\infty), \qquad x(\tau)=(1-\tau)x+\tau \hat x, \qquad \forall t\in [0,1], \ \tau\in [0,1].
\eeq
To this end, let $\tau\in(0,1]$, $\mathbb{j}=(j_1,\ldots,j_I)\in\mathcal{J}(x^\infty)$, and $x\in Y_\mathbb{j}(x^\infty)$ be arbitrarily chosen. It then follows from (\ref{cvl}) that
\begin{equation}\label{exafff111}
\phi_i(x)+\zeta_i(x)-\ell_{\psi_{i,j_i}}(x;x^\infty)\leq 0, \qquad \forall i\in\mathcal{I}.
\end{equation}
 By the convexity of $Y_\mathbb{j}(x^\infty)$, one has that $x(\tau)\in  Y_\mathbb{j}(x^\infty)$ and $x(t,\tau)\in  Y_\mathbb{j}(x^\infty)$ for all $t \in [0,1]$.} Also, by \eqref{slater-cond3}, \eqref{exafff111},  and the convexity of $\phi_i$ and $\zeta_i$, we have  
\begin{equation}\label{exahhh000}
\begin{array}{ll}
\hspace{-.15in}\phi_i(x(\tau))+\zeta_i(x(\tau))-\ell_{\psi_{i,j_i}}(x(\tau);x^\infty)\\
\hspace{-.15in} \le  (1-\tau) [\phi_i(x)+\zeta_i(x)-\ell_{\psi_{i,j_i}}(x;x^\infty)] + \tau [\phi_i(\hat x)+\zeta_i(\hat x)-\ell_{\psi_{i,j_i}}(\hat x;x^\infty)]
<0, \quad \forall i\in\mathcal{I}_=(x^\infty).
\end{array}
\end{equation}
Notice that $\phi_i(x^\infty)+\zeta_i(x^\infty)-\ell_{\psi_{i,j_i}}(x^\infty;x^\infty)=0$ for every $i\in\mathcal{I}_=(x^\infty)$. By this, the convexity of $\phi_i$ and $\zeta_i$, and a similar argument as for \eqref{exahhh000}, we obtain that for all $t \in [0,1]$ and $i\in\mathcal{I}_=(x^\infty)$,
$$
\begin{array}{ll}
\phi_i(x(t,\tau))+\zeta_i(x(t,\tau))-\ell_{\psi_{i,j_i}}(x(t,\tau);x^\infty)\\
\displaystyle \leq (1-t)\underbrace{[\phi_i(x^\infty)+\zeta_i(x^\infty)-\ell_{\psi_{i,j_i}}(x^\infty;x^\infty)]}_{=0}+ t[\phi_i(x(\tau))+\zeta_i(x(\tau))-\ell_{\psi_{i,j_i}}(x(\tau);x^\infty)]\\
= t[\phi_i(x(\tau))+\zeta_i(x(\tau))-\ell_{\psi_{i,j_i}}(x(\tau);x^\infty)].
\end{array}
$$
It follows from this, \eqref{hphi} and the convexity of $\phi_i$ that
\beq\label{uppbd-xtau1}
\ba{ll}
\hphi_i(x(t,\tau);x^\infty)+\zeta_i(x(t,\tau)) -\ell_{\psi_{i,j_i}}(x(t,\tau);x^\infty) \\
= \phi_i(x^\infty) + \nabla \phi_i(x^\infty)^T(x(t,\tau)-x^\infty)+L_i \|x(t,\tau)-x^\infty\|^2 /2+\zeta_i(x(t,\tau)) -\ell_{\psi_{i,j_i}}(x(t,\tau);x^\infty) \\
 \le  \phi_i(x(t,\tau))+L_it^2 \|x(\tau)-x^\infty\|^2/2 +\zeta_i(x(t,\tau)) -\ell_{\psi_{i,j_i}}(x(t,\tau);x^\infty)\\
\le t[\phi_i(x(\tau))+\zeta_i(x(\tau))-\ell_{\psi_{i,j_i}}(x(\tau);x^\infty)]
+L_it^2\|x(\tau)-x^\infty\|^2/2
\ea\eeq
for any $i\in\mathcal{I}_=(x^\infty)$ and $t \in [0,1]$.
By \eqref{exahhh000} and \eqref{uppbd-xtau1}, one can see that for $t>0$ sufficiently small,
\beq  \label{uppbd-xtau2}
\hphi_i(x(t,\tau);x^\infty)+\zeta_i(x(t,\tau)) -\ell_{\psi_{i,j_i}}(x(t,\tau);x^\infty) <0,\quad\forall i\in\mathcal{I}_=(x^\infty).
\eeq
 On the other hand, notice that $\phi_i(x^\infty)+\zeta_i(x^\infty)-\ell_{\psi_{i,j_i}}(x^\infty;x^\infty)<0$ for $i\not\in\mathcal{I}_=(x^\infty)$. By this, it is easy to observe
that \eqref{uppbd-xtau2} also holds for $i\not\in\mathcal{I}_=(x^\infty)$ and $t>0$ sufficiently small. Hence, there exists some $\hat t \in (0,1)$ such that the inequality
\eqref{uppbd-xtau3} holds.

{ In what follows, we show that
\beq\label{min_convex}
x^\infty \in \Argmin_{x\in Y_\mathbb{j}(x^\infty)} \ \underbrace{\hat{\phi}_0(x;x^\infty)+\zeta_0(x)-\ell_{\psi_{0,j_0}}(x;x^\infty)}_{\hat h_{0,j_0}(x)}, \qquad \forall j_0\in\mathcal{J}_0(x^\infty).
\eeq
}
Indeed, recall that $\{ x^k\}_{k\in \cK} \to x^\infty$. By \eqref{uppbd-xtau3} and the continuity of $\hphi_i(x(\hat t,\tau);\cdot)$ and $\ell_{\psi_{i,j_i}}(x(\hat t,\tau);\cdot)$,  one has that for sufficiently large $k\in\cK$,
\[
\hphi_i(x(\hat t,\tau);x^k)+\zeta_i(x(\hat t,\tau)) -\ell_{\psi_{i,j_i}}(x(\hat t,\tau);x^k) <0, \quad \forall i \in \cI.
\]
Replacing $x$ by $x(\hat t,\tau)$ in (\ref{exaaaa111}) and using this inequality,
 we can obtain that for every $j_0\in\mathcal{J}_0(x^\infty)$ and sufficiently large
$k\in \cK$,
\beq \label{bound-fval}
\phi_0(x^k)+\zeta_0(x^k)-\psi_0(x^k)-\eta_k
\leq \hat{\phi}_0(x(\hat t,\tau);x^k)+\zeta_0(x(\hat t,\tau))-\ell_{\psi_{0,j_0}}(x(\hat t,\tau);x^k).
\eeq
Taking limit on both sides of this inequality as $\cK \ni k \to\infty$ yields
\begin{equation}\label{lsls}
\phi_0(x^\infty)+\zeta_0(x^\infty)-\psi_0(x^\infty)
\leq \hat{\phi}_0(x(\hat t,\tau);x^\infty)+\zeta_0(x(\hat t,\tau))-\ell_{\psi_{0,j_0}}(x(\hat t,\tau);x^\infty)
\end{equation}
for any $j_0\in\mathcal{J}_0(x^\infty)$. By the convexity of $\hat{\phi}_0(\cdot;x^\infty)$ and $\zeta_0$, one has
\[
\ba{rcl}
\hat{\phi}_0(x(\hat t,\tau);x^\infty) &\leq& \hat t \hat{\phi}_0(x(\tau);x^\infty)+ (1-\hat t)\hat{\phi}_0(x^\infty;x^\infty), \\
\zeta_0 (x(\hat t,\tau)) &\leq & \hat t \zeta_0(x(\tau))+ (1-\hat t)\zeta_0(x^\infty).
\ea
\]
These, along with \eqref{exahhh333}, (\ref{lsls}), $\phi_0(x^\infty)=\hphi_{0}(x^\infty;x^\infty)$ and $\psi_0(x^\infty)=\ell_{\psi_{0,j_0}}(x^\infty;x^\infty)$ and the linearity of $\ell_{\psi_{0,j_0}}$, imply that for any fixed $\tau \in (0,1]$,
$$\hat{\phi}_0(x^\infty;x^\infty)+\zeta_0(x^\infty)-\ell_{\psi_{0,j_0}}(x^\infty;x^\infty) \leq \hat{\phi}_0(x(\tau);x^\infty)+\zeta_0(x(\tau))-\ell_{\psi_{0,j_0}}(x(\tau);x^\infty)$$
for any $j_0\in\mathcal{J}_0(x^\infty)$. Taking limit on both sides of this inequality by letting $\tau\downarrow 0$ gives
$$\hat{\phi}_0(x^\infty;x^\infty)+\zeta_0(x^\infty)-\ell_{\psi_{0,j_0}}(x^\infty;x^\infty) \leq \hat{\phi}_0(x;x^\infty)+\zeta_0(x)-\ell_{\psi_{0,j_0}}(x;x^\infty)$$
for any $j_0\in\mathcal{J}_0(x^\infty)$. Recall that $x$ is an arbitrary point in $Y_\mathbb{j}(x^\infty)$. It then follows from the above inequality that { \eqref{min_convex} holds.

We are now ready to complete the proof of statement (iii). Indeed, it follows from \eqref{min_convex} } that  $\hat h_{0,j_0}'(x^\infty;d) \ge 0$ for all $j_0\in\mathcal{J}_0(x^\infty)$ and $d\in\cT_{Y_\mathbb{j}(x^\infty)}(x^\infty)$. Notice from \eqref{Danskin} that for each $d\in\cT_{Y_\mathbb{j}(x^\infty)}(x^\infty)$ there
exists some $j_0\in\mathcal{J}_0(x^\infty)$ such that $F'(x^\infty; d) = \hat h_{0,j_0}'(x^\infty;d)$. It thus follows that $F'(x^\infty; d) \ge 0$ for all $d\in\cT_{Y_\mathbb{j}(x^\infty)}(x^\infty)$. By this, the arbitrariness of $\mathbb{j}\in\mathcal{J}(x^\infty)$, and the assumption that the PSCQ holds for $\Omega$ at $x^\infty$, we conclude from Theorem \ref{suff-cond-B} that $x^\infty$ is a B-stationary point of (\ref{e1}).

(iv) Suppose that $x^{\infty}\in\Omega$, $\mathcal{I}_=(x^\infty)\neq\emptyset$, $X$ is a polyhedral set, and moreover,  for every $i\in\mathcal{I}_=(x^{\infty})$ and $j_i\in\cJ_i(x^\infty)$, $\phi_{i}$ and $\psi_{i,j_i}$ are affine and $\zeta_i$ is piecewise affine on $X$. Recall from \eqref{hphi} that for every $i\in\cI$, $L_i=0$ if $\phi_i$ is affine. These together with \eqref{lg} and \eqref{hphi} imply that
\begin{equation}
\label{eq:affine}
\hat{\phi}_i(\tilde{x};\bx) = \phi_i(\tilde{x}), \quad \ell_{\psi_{i,j_i}}(\tilde{x};\bx) = \psi_{i,j_i}(\tilde{x}), \quad \forall \tilde{x},\bx\in X, \, i\in\cI_=(x^\infty), \, j_i\in\cJ_i(x^\infty).
\end{equation}
Let $\mathbb{j} = (j_1,\dots,j_I)\in\cJ(x^\infty)$ be arbitrarily chosen. For any $x\in Y_{\mathbb{j}}(x^\infty)$, we obtain from \eqref{cvl} that $x\in X$ and $\phi_i(x)+\zeta_i(x)-\ell_{\psi_{i,j_i}}(x;x^\infty) \leq 0$ for all $i\in\cI$. This together with $x^k,x^\infty\in X$ and \eqref{eq:affine} implies that
$$ \hat{\phi}_i(x;x^k) + \zeta_i(x) - \ell_{\psi_{i,j_i}}(x;x^k) \leq 0, \quad \forall i\in \cI_=(x^\infty), \ x\in Y_{\mathbb{j}}(x^\infty). $$
By this and a similar argument as in the proof {of \eqref{ee1}}, we obtain that for any $x\in Y_{\mathbb{j}}(x^\infty)$,
\beq\label{st_4}
\hat{\phi}_i(x^\infty+t (x-x^\infty);x^k)+\zeta_i(x^\infty+t (x-x^\infty))-\ell_{\psi_{i,j_i}}(x^\infty+t (x-x^\infty);x^k)\leq0,\quad\forall i\in\cI \nn
\eeq
holds for all $t$ sufficiently small and $k\in\cK$ sufficiently large. {Using this and a similar argument as for showing that $F'(x^\infty;d) \ge 0$ for all $d\in \mathcal{T}_X(x^\infty)$ in the proof of statement (ii), we have} that
$F'(x^\infty; d) \ge 0$ for all $d\in\cT_{Y_\mathbb{j}(x^\infty)}(x^\infty)$. By this, the arbitrariness of $\mathbb{j}\in\mathcal{J}(x^\infty)$, and the assumption that $X$ is a polyhedral set, $\phi_{i}$ is affine and $\zeta_i$ is piecewise affine on $X$ for every $i\in\mathcal{I}_=(x^\infty)$, we conclude from Theorem \ref{suff-cond-B} that $x^\infty$ is a B-stationary point of (\ref{e1}).
\end{proof}

\begin{remark}
\begin{itemize}
\item[i)] The assumptions in Theorem \ref{inexctct} are a natural generalization of the standard assumptions in the literature for classical penalty method for solving smooth constrained optimization problems (e.g., see~\cite[Theorem 17.2]{Nocedal}). In fact, they will be reduced to the standard assumptions when applied to the latter problems.
\item[ii)] For the case where the accumulation point $x^\infty$ of the solution sequence satisfies $x^\infty\in\Omega$ and $\cI_=(x^\infty)\neq\emptyset$, which is the most sophisticated case due to the presence of active DC constraints, we show in Theorem \ref{inexctct} (iii) that $x^\infty$ is a B-stationary point of \eqref{e1}, provided that the PSCQ in Definition \ref{slater} holds for $\Omega$ at $x^\infty$. In contrast with the convergence result presented in \cite[Proposition 9]{Pang}, our result does not require the additional assumption that the set $\{ j:  \psi_{i,j}(x^\infty) = \psi_i(x^\infty)\}$ is a singleton for $i=1,\ldots, I$. As a consequence, the feasible accumulation point $x^\infty$ of the solution sequence generated by our penalty method can be a B-stationary point of \eqref{e1} even when some of $\psi_i$'s are non-differentiable at $x^\infty$. In addition, as we have mentioned in Remark \ref{B-statpt}, the PSCQ condition used in our result is  generally weaker than that in \cite{Pang}.
\end{itemize}
\end{remark}

\setcounter{equation}{0}
\section{An augmented Lagrangian method for DC program (\ref{e1})}
\label{alm}


In this section we propose an augmented Lagrangian (AL) method for solving problem \eqref{e1} and analyze its convergence. We also provide an example to demonstrate its convergence. To this end,  we introduce an AL function for \eqref{e1} given by
\beq\label{al-fun}
\tF_{\rho}(x,\lambda)=F(x)+\frac{1}{2\rho}\sum\limits_{i=1}^I([\lambda_i+\rho(\phi_i(x)+\zeta_i(x)-\psi_i(x))]_+^2-\lambda_i^2),
\eeq
where $\rho>0$. For any $\bar{x}\in X$, $\lambda\in\Re^I$, $j_0\in \mathcal{J}_0$, and $\mathbb{j}\in\mathcal{J}$, we define $\tQ_\rho(\cdot;\bar{x},\lambda,j_0,\mathbb{j}): X\to \Re$ by
\beq\label{up-al-fun}
\tQ_\rho(x;\bar{x},\lambda,j_0,\mathbb{j})=\hphi_0(x;\bar{x}) + \zeta_0(x) - \ell_{\psi_{0,j_0}}(x; \bx) + \frac{1}{2\rho}\sum\limits_{i=1}^I([\lambda_i+\rho(\hphi_i(x;\bar{x}) + \zeta_i(x) - \ell_{\psi_{i,j_i}}(x; \bx))]_+^2-\lambda_i^2),
\eeq
where $\ell_{\psi_{i,j_i}}$ and $\hphi_i$ are defined in \eqref{lg} and \eqref{hphi}, respectively.
It is easy to see from (\ref{lg}) and (\ref{hphi}) that for every $j_0\in\cJ_0(\bar{x})$ and $\bj\in\cJ(\bar{x})$, $\tF_{\rho}(x,\lambda) \leq \tQ_\rho(x;\bar{x},\lambda,j_0,\mathbb{j})$ for all $x\in X$ and $\tQ_\rho(\cdot;\bar{x},\lambda,j_0,\mathbb{j})$ is strongly convex on $X$ with modulus $L_0>0$.

We now propose an AL method for solving problem \eqref{e1} in which a sequence of AL subproblems are approximately solved.  The details of the AL method are presented as follows.

\begin{algorithm}
\label{aug-lag}
\normalfont
\mbox{}
\begin{itemize}
\item[0.] Input $\epsilon>0$, $\rho_0>0$, $\alpha>0$, $\sigma>1$, $\lambda^0\in \Re^I_+$, and a sequence $\{\eta_k\}\subset \Re_+$ such that $\eta_k\to0$. Set $k\leftarrow 0$.
\item[1.] Find an approximate solution $x^k$ of the AL subproblem
\begin{equation}\label{al-pro1}
\min\limits_{x\in X}\tF_{\rho_k}(x,\lambda^k)
\end{equation}
such that $x^k\in X$ and
\begin{equation}\label{fun-value}
\tF_{\rho_k}(x^{k},\lambda^k)\leq \tQ_{\rho_k}(x;x^k,\lambda^k,j_0,\mathbb{j})+\eta_k,\qquad\forall x\in X
\end{equation}
for every $j_0\in\mathcal{J}_{0,\epsilon}(x^{k})$ and $\mathbb{j}\in\mathcal{J}_{\epsilon}(x^{k})$, where $\cJ_{0,\epsilon}(x^{k})$ and $\cJ_{\epsilon}(x^{k})$ are defined in (\ref{J0-index}) and (\ref{J-index}), respectively.
\item[2.] Update $\lambda^{k+1}=(\lambda^{k+1}_1,\ldots,\lambda^{k+1}_I)^T$ by
$$\lambda_i^{k+1} =\left[\lambda^k_i+\rho_k\left(\phi_i(x^{k})+\zeta_i(x^k)-\psi_{i}(x^{k})\right)\right]_+,\quad\forall i=1,\ldots,I.$$
\item[3.] Set $\rho_{k+1} = \max\{\sigma \rho_k,\|\lambda^{k+1}\|^{1+\alpha}\}$.
\item[4.] Set $k\leftarrow k+1$, and go to Step 1.
\end{itemize}
{\bf End.}
\end{algorithm}

\begin{remark}
\bi
\item[(i)] The approximate solution $x^k$ of \eqref{al-pro1} satisfying $x^k \in X$ and \eqref{fun-value} can be found by Algorithm \ref{I_SCA} proposed in Section \ref{penalty-subprob}.

\item[(ii)] The update scheme on penalty parameters is adopted from \cite{Lu12}, which differs from the one for the classical AL method in that the magnitude of the penalty parameters in our method outgrows that of Lagrangian multipliers.
\ei
\end{remark}

We next establish a  convergence result for Algorithm \ref{aug-lag}.

\begin{theorem}\label{alm-them}
Let $\{x^k\}$ be generated by Algorithm \ref{aug-lag}. Assume that $\{x^k\}_{k\in \mathcal{K}}$ converges to $x^\infty$ for some subsequence $\cK$. Then the following statements hold.
\bi
\item[(i)] $x^\infty$ is a D-stationary point of the problem
$$\min\limits_{x\in X}\sum\limits_{i=1}^I[\phi_i(x)+\zeta_i(x)-\psi_{i}(x)]_+^2.$$
\item[(ii)] If $x^{\infty}\in\Omega$ and $\mathcal{I}_=(x^\infty)=\emptyset$, then $x^\infty$ is a D-stationary point of $\min\limits_{x\in X} F(x)$.
\item[(iii)] If $x^{\infty}\in\Omega$, $\mathcal{I}_=(x^\infty)\neq\emptyset$, and moreover, the PSCQ holds for $\Omega$ at $x^\infty$, then $x^\infty$ is a B-stationary point of problem (\ref{e1}).
\item[(iv)] Suppose that $\alpha>1$ in Algorithm \ref{aug-lag}. If $x^{\infty}\in\Omega$, $\mathcal{I}_=(x^\infty)\neq\emptyset$, $X$ is a polyhedral set, and moreover, for every $i\in\mathcal{I}_=(x^{\infty})$ and $j_i\in\cJ_i(x^\infty)$, $\phi_{i}$ and $\psi_{i,j_i}$ are affine and $\zeta_i$ is piecewise affine on $X$, then $x^\infty$ is a B-stationary point of problem (\ref{e1}).
\ei
\end{theorem}

\begin{proof}
Since $\{x^k\}_{k\in \cK}\to x^\infty$, one has  $\mathcal{J}_0(x^\infty)\subseteq \mathcal{J}_{0,\epsilon}(x^k)$ and $\mathcal{J}(x^\infty)\subseteq \mathcal{J}_\epsilon(x^k)$ for sufficiently large $k\in \cK$. It thus follows from \eqref{al-fun}, \eqref{up-al-fun} and \eqref{fun-value} that for $k\in\cK$ sufficiently large, one has
\begin{equation}\label{fun-value-2}
\begin{array}{ll}
\displaystyle{\phi_0(x^k)+\zeta_0(x^k)-\psi_0(x^k)+\frac{1}{2\rho_k}\sum\limits_{i=1}^I\left[\lambda_i^k+\rho_k\left(\phi_i(x^k)+\zeta_i(x^k)-\psi_i(x^k)\right)\right]_+^2}\\
\displaystyle\leq \hat{\phi}_0(x;x^k) +\zeta_0(x)-\ell_{\psi_{0,j_0}}(x;x^k)+\frac{1}{2\rho_k}\sum\limits_{i=1}^I\left[\lambda^k_i+\rho_k\left(\hat{\phi}_i(x;x^k)+\zeta_i(x)-\ell_{\psi_{i,j_i}}(x;x^k)\right)\right]_+^2+\eta_k,
\end{array}
\end{equation}
for all $j_0\in\mathcal{J}_0(x^\infty)$, $\mathbb{j}=(j_1,\ldots,j_I)\in\mathcal{J}(x^\infty)$ and $x\in X$. In addition, one can observe from Step 3 of Algorithm \ref{aug-lag} that $\{\rho_k\} \to \infty$ and $\{\lambda^k/\rho_k\} \to 0$.

(i) Dividing both sides of \eqref{fun-value-2} by $\rho_k$,
taking limits as $\cK\ni k \to \infty$, and using $\{x^k\}_{k\in \cK}\to x^\infty$, $\{\rho_k\} \to \infty$ and $\{\lambda^k/\rho_k\} \to 0$, we have
\[
\sum\limits_{i=1}^I[\phi_i(x^\infty)+\zeta_i(x^\infty)-\psi_i(x^\infty)]_+^2 \le \sum\limits_{i=1}^I[\hat{\phi}_i(x;x^\infty)+\zeta_i(x)-\ell_{\psi_{i,j_i}}(x;x^\infty)]_+^2,\ \ \forall x\in X
\]
for any $\mathbb{j}=(j_1,\ldots,j_I) \in \mathcal{J}(x^\infty)$. The rest of the proof of this statement follows from this inequality, and
{the similar arguments as the ones that are from \eqref{constr-opt} till the end of the proof of Theorem \ref{inexctct} (i).}

(ii)  Assume that $x^{\infty}\in\Omega$ and $\mathcal{I}_=(x^\infty)=\emptyset$. Let $x\in X$ be arbitrarily chosen. By a similar argument as in the proof {of \eqref{ee1}}, one can show that there exists some $\delta<0$ such that for any $i\in\mathcal{I}$, $\phi_i(x^k)+\zeta_i(x^k)-\psi_{i}(x^k)<\delta$ and
\[
\hat{\phi}_i(x^\infty+t (x-x^\infty);x^k)+\zeta_i(x^\infty+t (x-x^\infty))-\ell_{\psi_{i,j_i}}(x^\infty+t (x-x^\infty);x^k)<\delta,\ \forall  j_i\in\mathcal{J}_i(x^\infty)
\]
hold for all $k\in \cK$ sufficiently large and $t>0$ sufficiently small. By these two relations and the fact $\{\rho_k\}\to \infty$ and $\{\lambda^k/\rho_k\}\to 0$, one can obtain that for all $k\in \cK$ sufficiently large and $t>0$ sufficiently small,
$[\lambda_i^k+\rho_k(\phi_i(x^k)+\zeta_i(x^k)-\psi_i(x^k))]_+=0$ and
\[
\left[\lambda^k_i+\rho_k\left(\hat{\phi}_i(x^\infty+t (x-x^\infty);x^k)+\zeta_i(x^\infty+t (x-x^\infty))-\ell_{\psi_{i,j_i}}(x^\infty+t (x-x^\infty);x^k)\right)\right]_{+}=0.
\]
Using these two relations and replacing $x$ by $x^\infty+t (x-x^\infty)$ in (\ref{fun-value-2}), we have that  for every $x\in X$ and $j_0\in\mathcal{J}_0(x^\infty)$,
\[\phi_0(x^k)+\zeta_0(x^k)-\psi_0(x^k)\leq \hat{\phi}_0(x^\infty+t (x-x^\infty);x^k)+\zeta_0(x^\infty+t (x-x^\infty))-\ell_{\psi_{0,j_0}}(x^\infty+t (x-x^\infty);x^k)+\eta_k\]
for all $k\in \cK$ sufficiently large and $t>0$ sufficiently small.  The rest of the proof of this statement follows from this inequality, and
{the similar arguments as the ones that are from \eqref{ee1} till the end of the proof of Theorem \ref{inexctct} (ii).}

(iii) Assume that $x^{\infty}\in\Omega$, $\mathcal{I}_=(x^\infty)\neq\emptyset$, and moreover,  the PSCQ holds for $\Omega$ at $x^\infty$. Let $\mathbb{j}=(j_1,\ldots,j_I)\in\mathcal{J}(x^\infty)$ and $x\in Y_\mathbb{j}(x^\infty)$ be arbitrarily chosen, and let $x(t,\tau)$ be defined in \eqref{exahhh333} for all $t,\tau\in [0,1]$.
By a similar argument as in the proof of Theorem \ref{inexctct} (iii),
one can show that for any fixed $\tau\in(0,1]$, there exist some $\hat{t}\in(0,1)$ and $\delta<0$ that are dependent on $\tau$ such that $\hphi_i(x(\hat{t},\tau);x^k)+\zeta_i(x(\hat{t},\tau)) -\ell_{\psi_{i,j_i}}(x(\hat{t},\tau);x^k) <\delta$
for all $i \in \mathcal{I}$ and $k\in\mathcal{K}$ sufficiently large.  It then follows from this relation, $\{\rho_k\}\to\infty$ and $\{\lambda^k/\rho_k\}\to 0$ that for all $i \in \mathcal{I}$ and $k\in\mathcal{K}$ sufficiently large,
$[\lambda_i^k+\rho_k(\hphi_i(x(\hat{t},\tau);x^k)+\zeta_i(x(\hat{t},\tau)) -\ell_{\psi_{i,j_i}}(x(\hat{t},\tau);x^k) )]_+=0$.
Replacing $x$ by $x(\hat{t},\tau)$ in \eqref{fun-value-2} and using this relation, we see that \eqref{bound-fval} holds for every $j_0\in\mathcal{J}_0(x^\infty)$ and sufficiently large
$k\in \cK$. The rest of the proof of this statement follows from \eqref{bound-fval}, and
{the similar arguments as the ones that are from \eqref{bound-fval} till the end of the proof of Theorem \ref{inexctct} (iii).}

(iv) From $\alpha>1$ and Step 3 of Algorithm \ref{alm}, we observe that $\|\lambda^k\|^2/\rho_k\rightarrow 0$ as $k\rightarrow \infty$. By a similar argument as in the proof of Theorem \ref{inexctct} (iv), one can show that
\[\left[\lambda^k_i+\rho_k\left(\hat{\phi}_i(x^\infty+t (x-x^\infty);x^k)+\zeta_i(x^\infty+t (x-x^\infty))-\ell_{\psi_{i,j_i}}(x^\infty+t (x-x^\infty);x^k)\right)\right]_+\leq \lambda^k_i\]
for any $x\in Y_\mathbb{j}(x^\infty)$, $i\in\mathcal{I}$, $j_i\in \mathcal{J}_i(x^\infty)$, $k\in \cK$ sufficiently large and $t>0$ sufficiently small. Replacing $x$ by $x^\infty+t (x-x^\infty)$ in (\ref{fun-value-2}), we have that  for every $x\in Y_\mathbb{j}(x^\infty)$ and $j_0\in\mathcal{J}_0(x^\infty)$,
\[\ba{l}\phi_0(x^k)+\zeta_0(x^k)-\psi_0(x^k)\\
\leq\displaystyle \hat{\phi}_0(x^\infty+t (x-x^\infty);x^k)+\zeta_0(x^\infty+t (x-x^\infty))-\ell_{\psi_{0,j_0}}(x^\infty+t (x-x^\infty);x^k)+\frac{\|\lambda^k\|^2}{2\rho_k}+\eta_k
\ea\]
for all $k\in \cK$ sufficiently large and $t>0$ sufficiently small. Taking limits on both sides of this inequality by letting $\cK\ni k \to \infty$, and using $\{x^k\}_{k\in \cK}\to x^\infty$, $\rho_k \to \infty$ and $\|\lambda^k\|^2/\rho_k \to 0$ as $k\rightarrow\infty$, it follows that for every $x\in Y_\mathbb{j}(x^\infty)$, $j_0\in\mathcal{J}_0(x^\infty)$, and $t>0$ sufficiently small,
\[\phi_0(x^\infty)+\zeta_0(x^\infty)-\psi_0(x^\infty)\leq \hat{\phi}_0(x^\infty+t (x-x^\infty);x^\infty)+\zeta_0(x^\infty+t (x-x^\infty))-\ell_{\psi_{0,j_0}}(x^\infty+t (x-x^\infty);x^\infty).\]
{ By this and the similar arguments as the ones that are from \eqref{lim-ineq} till the end of the proof of Theorem \ref{inexctct} (ii), one can obtain
 that $F'(x^\infty; d) \ge 0$ for all $d\in\cT_{Y_\mathbb{j}(x^\infty)}(x^\infty)$. By this, the arbitrariness of $\mathbb{j}\in\mathcal{J}(x^\infty)$, and the assumption that $X$ is a polyhedral set, $\phi_{i}$ is affine and $\zeta_i$ is piecewise affine on $X$ for every $i\in\mathcal{I}_=(x^\infty)$, we conclude from Theorem \ref{suff-cond-B} that $x^\infty$ is a B-stationary point of (\ref{e1}).}
\end{proof}

The above theorem establishes the convergence of $\{x^k\}$. Nevertheless,  the convergence of $\{\lambda^k\}$ remains unknown. Even if $\{\lambda^k\}$ converges (subsequentially) to some $\lambda^\infty$, it appears impossible to satisfy the KKT conditions \eqref{kkt-1} and \eqref{kkt-2}. In the next theorem, we construct a set of auxiliary Lagrangian multiplier sequences and show that their accumulation points together with the accumulation points of $\{x^k\}$ satisfy the KKT conditions \eqref{kkt-1} and \eqref{kkt-2}.

\begin{theorem}\label{alm-them-mtp}
Let $\{x^k\}$ and $\{\lambda^k\}$ be generated by Algorithm \ref{aug-lag}. Suppose that $\{x^k\}_{k\in \mathcal{K}}$ converges to $x^\infty$ for some subsequence $\cK$.  Assume  that $x^{\infty}\in\Omega$, $\mathcal{I}_=(x^\infty)\neq\emptyset$, {$\zeta_0(x^\infty)<\infty$}, and the PSCQ holds for $\Omega$ at $x^\infty$.  For any $j_0\in\cJ_0(x^\infty)$ and $\bj=(j_1,j_2,\ldots,j_I)\in\cJ(x^\infty)$, suppose that $x^{k,j_0,\bj}\in X$ satisfies
\beq\label{dist-subdif}
{\rm dist}\left(0, \partial [\tQ_{\rho_k}(x;x^{k},\lambda^k,j_0,\mathbb{j})+\iota_X(x)]\big|_{x=x^{k,j_0,\mathbb{j}}}\right)\le\gamma_k
\eeq
with $\{\gamma_k\} \to 0$.
Let $\lambda^{k,j_0,\bj}=(\lambda_1^{k,j_0,\bj},\lambda_2^{k,j_0,\bj},\ldots,\lambda_I^{k,j_0,\bj})$,
where
\beq\label{lam-def}
\lambda_i^{k,j_0,\bj}=\left[\lambda^k_i+\rho_k\left(\hat{\phi}_i(x^{k,j_0,\bj};x^k)+\zeta_i(x^{k,j_0,\bj})-\ell_{\psi_{i,j_i}}(x^{k,j_0,\bj};x^k)\right)\right]_+.
\eeq
Then the following statements hold.
\bi
\item[(i)]  $\{x^{k,j_0,\bj}\}_{k\in \cK}$ converges to $x^\infty$.
\item[(ii)] $\{\lambda^{k,j_0,\bj}\}_{k\in \cK}$ is bounded. Moreover, every accumulation point $\lambda^{\infty,j_0,\bj}$ of  $\{\lambda^{k,j_0,\bj}\}_{k\in \cK}$ satisfies that
\[
\ba{l}
\lambda_i^{\infty,j_0,\bj}\geq0,\quad \lambda_i^{\infty,j_0,\bj}\left[\phi_i(x^\infty)+\zeta_i(x^\infty)-\psi_{i,j_i}(x^\infty)\right]=0,\quad\forall i\in\cI, \\ 
0\in \nabla\phi_0(x^\infty)+\partial \zeta_0(x^\infty)-\nabla\psi_{0,j_0}(x^\infty)\\
\quad\ \ +\sum\limits_{i=1}^I\lambda_i^{\infty,j_0,\bj}[\nabla\phi_i(x^\infty)+\partial\zeta_i(x^\infty)-\nabla\psi_{i,j_i}(x^\infty)]+\mathcal{N}_X(x^\infty).
\ea
\]
\ei
\end{theorem}

\begin{proof}{ (i) In order to prove statement (i), we first show that $\{x^{k,j_0,\mathbb{j}}\}_{k\in \mathcal{K}}$ is bounded. To this end, let} $j_0\in\cJ_0(x^\infty)$ and $\bj=(j_1,j_2,\ldots,j_I)\in\cJ(x^\infty)$ be arbitrarily chosen. One can observe from (\ref{hphi}) and (\ref{up-al-fun}) that $\tQ_{\rho_k}(\cdot;x^k,\lambda^k,j_0,\mathbb{j})$ is strongly convex on $X$ with modulus $L_0>0$. Since $x^{k,j_0,\bj}$ satisfies (\ref{dist-subdif}), there exists some $s\in \partial [\tQ_{\rho_k}(x;x^k,\lambda^k,j_0,\mathbb{j})+\iota_X(x)]\big|_{x=x^{k,j_0,\mathbb{j}}}$ such that $\|s\| \le \gamma_k$. This along with the strong convexity of $\tQ_{\rho_k}(\cdot;x^k,\lambda^k,j_0,\mathbb{j})$ yields
\[\ba{lll}
\tQ_{\rho_k}(x;x^k,\lambda^k,j_0,\mathbb{j})&\ge& \min\limits_z \tQ_{\rho_k}(x^{k,j_0,\bj};x^k,\lambda^k,j_0,\mathbb{j}) +s^T(z-x^{k,j_0,\bj}) + L_0\|z-x^{k,j_0,\bj}\|^2/2\\
& =& \tQ_{\rho_k}(x^{k,j_0,\bj};x^k,\lambda^k,j_0,\mathbb{j}) -  \|s\|^2/(2L_0)\\
 &\ge& \tQ_{\rho_k}(x^{k,j_0,\bj};x^k,\lambda^k,j_0,\mathbb{j}) - \gamma_k^2/(2L_0),\quad \forall x\in X,
\ea
\]
which gives
\beq\label{app-opt}
\tQ_{\rho_k}(x^{k,j_0,\bj};x^k,\lambda^k,j_0,\mathbb{j})\le \tQ_{\rho_k}(x;x^k,\lambda^k,j_0,\mathbb{j})+\frac{\gamma_k^2}{2L_0},\qquad \forall x\in X.
\eeq
Note that $x^{\infty}\in\Omega $, $\mathcal{I}_=(x^\infty)\neq\emptyset$, and moreover, the PSCQ holds for $\Omega$ at $x^\infty$. {By a similar argument as for proving \eqref{uppbd-xtau3},}   one can show that there exist $\bar{t}\in(0,1)$ and $\hat{x}\in Y_{\bj}(x^\infty)$ satisfying \eqref{slater-cond3}
 such that $\hphi_i(\hat{x}(t);x^\infty)+\zeta_i(\hat{x}(t)) -\ell_{\psi_{i,j_i}}(\hat{x}(t);x^\infty) <0$ for any $i \in \mathcal{I}$ and $t\in(0,\bar{t})$,
where $\hat{x}(t)=x^\infty+t(\hat{x}-x^\infty)$. Let us fix $t\in(0,\bar t)$ arbitrarily. 
Observe that $\hphi_i(\cdot;\cdot)$ and $\ell_{\psi_{i,j_i}}(\cdot;\cdot)$ are continuous on $X\times X$. It then follows from the above relation, $\{\rho_k\}\to\infty$ and $\{\lambda^k/\rho_k\}\to 0$ that for any $k\in\mathcal{K}$ sufficiently large and $i \in \mathcal{I}$, $[\lambda_i^k+\rho_k(\hphi_i(\hat{x}(t);x^k)+\zeta_i(\hat{x}(t)) -\ell_{\psi_{i,j_i}}(\hat{x}(t);x^k))]_+=0$.
Replacing $x$ by $\hat{x}(t)$ in (\ref{app-opt}) and using this equality, we can  obtain from (\ref{up-al-fun}) that for all $k\in \mathcal{K}$ sufficiently large, one has
\beq\label{fun-rela}
\ba{ll}
\hat{\phi}_0(x^{k,j_0,\mathbb{j}};x^k)+\zeta_0(x^{k,j_0,\mathbb{j}})-\ell_{\psi_{0,j_0}}(x^{k,j_0,\mathbb{j}};x^k) \\
\displaystyle+ \frac{1}{2\rho_k}\sum\limits_{i=1}^I
\left[\lambda^k_i+\rho_k\left(\hphi_i(x^{k,j_0,\mathbb{j}};x^k) + \zeta_i(x^{k,j_0,\mathbb{j}}) - \ell_{\psi_{i,j_i}}(x^{k,j_0,\mathbb{j}}; x^k)\right)\right]_+^2\\
\displaystyle= \tQ_{\rho_k}(x^{k,j_0,\bj};x^k,\lambda^k,j_0,\mathbb{j})+\frac{1}{2\rho_k}\sum\limits_{i=1}^I(\lambda_i^k)^2\\
\leq\displaystyle \hat{\phi}_0(\hat{x}(t);x^k)+\zeta_0(\hat{x}(t))-\ell_{\psi_{0,j_0}}(\hat{x}(t);x^k)+\frac{\gamma_k^2}{2L_0},
\ea
\eeq
which implies that
\beq\label{fun-rela1}
\hat{\phi}_0(x^{k,j_0,\mathbb{j}};x^k)+\zeta_0(x^{k,j_0,\mathbb{j}})-\ell_{\psi_{0,j_0}}(x^{k,j_0,\mathbb{j}};x^k) \leq \hat{\phi}_0(\hat{x}(t);x^k)+\zeta_0(\hat{x}(t))-\ell_{\psi_{0,j_0}}(\hat{x}(t);x^k)+\frac{\gamma_k^2}{2L_0}.
\eeq
Claim that $\{x^{k,j_0,\mathbb{j}}\}_{k\in \mathcal{K}}$ is bounded. Suppose for contradiction that $\{x^{k,j_0,\mathbb{j}}\}_{k\in \mathcal{K}}$ is unbounded. By passing to a subsequence if necessary,  we  assume that $\lim\limits_{\cK\ni k\to\infty}\|x^{k,j_0,\mathbb{j}}\|=\infty$.  By \eqref{e2}, \eqref{hphi},  \eqref{fun-rela1} and the convexity of $\zeta_0$, we have
\[\ba{ll}
 \displaystyle F(x^k)-\|\nabla\phi_0(x^k)+v_0^k-\nabla\psi_{0,j_0}(x^k)\|\|x^{k,j_0,\bj}-x^k\|+\frac{L_0}{2}\|x^{k,j_0,\mathbb{j}}-x^k\|^2\\
 \displaystyle\leq \hat{\phi}_0(x^{k,j_0,\mathbb{j}};x^k) +\zeta_0(x^{k,j_0,\mathbb{j}})-\ell_{\psi_{0,j_0}}(x^{k,j_0,\mathbb{j}};x^k)
\leq\hat{\phi}_0(\hat{x}(t);x^k)+\zeta_0(\hat{x}(t))-\ell_{\psi_{0,j_0}}(\hat{x}(t);x^k)+\frac{\gamma_k^2}{2L_0},
\ea
\]
where $v_0^k\in\partial\zeta_0(x^k)$. Since $\{x^k\}_{k\in \mathcal{K}}\to x^\infty$, it follows from \cite[Theorems 23.4 and 24.5]{Rockafellar} that $\cup_{k\in \cK}\partial\zeta_0(x^k)$ is bounded and so is $\{v_0^k\}_{k\in \cK}$. Using these,
$\{\gamma_k\}\to0$, $\lim\limits_{\cK\ni k\to\infty}\|x^{k,j_0,\mathbb{j}}\|=\infty$, and taking limit as $\cK\ni k\to\infty$ on both sides of the last inequality, we obtain  $\infty\leq \hat{\phi}_0(\hat{x}(t);x^\infty)+\zeta_0(\hat{x}(t))-\ell_{\psi_{0,j_0}}(\hat{x}(t);x^\infty)$, which clearly cannot hold. Hence, $\{x^{k,j_0,\mathbb{j}}\}_{k\in \mathcal{K}}$ is bounded.

{We are now ready to complete the proof of statement (i). Indeed,} since $\{x^{k,j_0,\mathbb{j}}\}_{k\in \mathcal{K}}$ is bounded, it suffices to show that each convergent subsequence of $\{x^{k,j_0,\mathbb{j}}\}_{k\in \mathcal{K}}$ converges to $x^\infty$. By passing to a subsequence if necessary, we can assume that $\{x^{k,j_0,\bj}\}_{k\in \mathcal{K}}\to x^{\infty, j_0,\bj}$. Using this and taking limit on both sides of (\ref{fun-rela1}) as $\cK\ni k\to\infty$, we have
\[
\hat{\phi}_0(x^{\infty,j_0,\mathbb{j}};x^\infty)+\zeta_0(x^{\infty,j_0,\mathbb{j}})-\ell_{\psi_{0,j_0}}(x^{\infty,j_0,\mathbb{j}};x^\infty)\leq \hat{\phi}_0(\hat{x}(t);x^\infty)+\zeta_0(\hat{x}(t))-\ell_{\psi_{0,j_0}}(\hat{x}(t);x^\infty)
\]
for any $t \in (0,\bar t)$. Taking limit on both sides of this inequality as $t\downarrow 0$ gives
\beq\label{opt1}
\hat{\phi}_0(x^{\infty,j_0,\mathbb{j}};x^\infty)+\zeta_0(x^{\infty,j_0,\mathbb{j}})-\ell_{\psi_{0,j_0}}(x^{\infty,j_0,\mathbb{j}};x^\infty) \leq \hat{\phi}_0(x^\infty;x^\infty)+\zeta_0(x^\infty)-\ell_{\psi_{0,j_0}}(x^\infty;x^\infty).
\eeq
Recall that $\{\rho_k\} \to \infty$, $\{\lambda^k/\rho_k\} \to 0$ and $\{\gamma_k\} \to 0$. Using these,  dividing both sides of \eqref{fun-rela} by $\rho_k$, and taking limit as $\cK\ni k\to\infty$, we obtain that $\sum\limits_{i=1}^I[\hphi_i(x^{\infty,j_0,\bj};x^\infty)+\zeta_i(x^{\infty,j_0,\bj}) -\ell_{\psi_{i,j_i}}(x^{\infty,j_0,\bj};x^\infty)]^2_+\leq0$,
which together with $\phi_i(x^{\infty,j_0,\bj})\le \hphi_i(x^{\infty,j_0,\bj};x^\infty)$ implies that for any $i\in \cI$, $\phi_i(x^{\infty,j_0,\bj})+\zeta_i(x^{\infty,j_0,\bj}) -\ell_{\psi_{i,j_i}}(x^{\infty,j_0,\bj};x^\infty)\leq0$.
Hence, $x^{\infty,j_0,\bj}\in Y_\bj(x^\infty)$. Recall that $x^\infty$ is a B-stationary point of (\ref{e1}) and the PSCQ holds for $\Omega$ at $x^\infty$. It then follows from Theorem \ref{suff-cond-B} that $F'(x^\infty; d) \ge 0$ for all $d\in \cT_{Y_\bj(x^\infty)}(x^\infty)$, which along with $j_0\in \cJ_0(x^\infty)$ implies that
\[
\nabla \phi_0(x^\infty)^T d + \zeta'_0(x^\infty;d) - \nabla \psi_{0,j_0}(x^\infty)^T d \ \ge \ F'(x^\infty; d) \ \ge \ 0, \quad \forall d\in \cT_{Y_\bj(x^\infty)}(x^\infty).
\]
Hence, by \eqref{hphi}, one has
\beq\label{opt3}
x^\infty=\argmin_{x\in Y_\bj(x^\infty)} \hat{\phi}_0(x;x^\infty)+\zeta_0(x)-\ell_{\psi_{0,j_0}}(x;x^\infty).
\eeq
Notice that the objective of \eqref{opt3} is strongly convex. Hence, $x^\infty$ is the unique optimal solution of \eqref{opt3}. In addition, observe from \eqref{opt1} and
$x^{\infty,j_0,\bj}\in Y_\bj(x^\infty)$ that $x^{\infty,j_0,\bj}$ is also an optimal solution of \eqref{opt3}. It then follows that $x^{\infty,j_0,\mathbb{j}}=x^\infty$. Therefore, statement (i) holds as desired.

(ii)  Let $j_0\in\cJ_0(x^\infty)$ and $\bj=(j_1,j_2,\ldots,j_I)\in\cJ(x^\infty)$ be chosen arbitrarily. It follows from (\ref{up-al-fun}) and (\ref{dist-subdif}) that there exist $v^{k,j_0,\bj}_i\in \partial\zeta_i(x^{k,j_0,\bj})$ for $i\in\{0,1,\ldots,I\}$  and $w^{k,j_0,\bj}\in \mathcal{N}_X(x^{k,j_0,\bj})$ such that
\begin{equation}\label{lg2222}
\begin{array}{ll}
\left\|\nabla \phi_0(x^k) + L_0(x^{k,j_0,\bj}-x^k)+v^{k,j_0,\bj}_0-\nabla \psi_{0,j_0}(x^k)\right.\\
\displaystyle\left.
+\sum\limits_{i=1}^I\lambda_i^{k,j_0,\bj}[\nabla \phi_i(x^k) + L_i(x^{k,j_0,\bj}-x^k)+v^{k,j_0,\bj}_i-\nabla \psi_{i,j_i}(x^k)]+w^{k,j_0,\bj}
\right\| \le \gamma_k,
\end{array}
\end{equation}
where $\lambda_i^{k,j_0,\bj}$ is defined in \eqref{lam-def}.

{In order to prove statement (ii), we first show that $\{\lambda^{k,j_0,\bj}\}_{k\in \cK}$ is bounded.} Suppose for contradiction that it is unbounded.  By passing to a subsequence if necessary, we can assume that $\{\|\lambda^{k,j_0,\bj}\|\}_{k\in \cK}\to \infty$.  Denote $\bar{\lambda}^{k,j_0,\bj}=\lambda^{k,j_0,\bj}/\|\lambda^{k,j_0,\bj}\|$ and $\bar{w}^{k,j_0,\bj}=w^{k,j_0,\bj}/\|\lambda^{k,j_0,\bj}\|$. Then $\|\bar{\lambda}^{k,j_0,\bj}\|=1$ for all $k\in \cK$. Recall that $\{x^k\}_{k\in \cK}\to x^\infty$ and $\{x^{k,j_0,\bj}\}_{k\in \cK}\to x^\infty$. It then follows from \cite[Theorems 23.4 and 24.5]{Rockafellar} that $\cup_{k\in \cK}\partial\zeta_i(x^{k,j_0,\bj})$ is bounded and so is $\{v_i^{k,j_0,\bj}\}_{k\in \cK}$ for all $i\in\{0,1,\ldots,I\}$. In addition, notice that $\{\nabla \phi_i(x^{k})\}_{k\in \cK}\to \nabla \phi_{i}(x^\infty)$, $\{\nabla \psi_{i,j_i}(x^{k})\}_{k\in \cK}\to \nabla \psi_{i,j_i}(x^\infty)$ and $\{\gamma_k\}\to0$. In view of these and  (\ref{lg2222}), one can observe that $\{\bar{w}^{k,j_0,\bj}\}$ is bounded. By passing to a subsequence if necessary, we can assume that $\{\bar{\lambda}^{k,j_0,\bj}\}_{k\in \cK}\to \bar{\lambda}^{\infty,j_0,\bj}$, $\{v_i^{k,j_0,\bj}\}_{k\in \cK}\to v_i^{\infty,j_0,\bj}$ for any $i\in\cI$ and $\{\bar{w}^{k,j_0,\bj}\}_{k\in \cK}\to \bar{w}^{\infty,j_0,\bj}$. Clearly, $\bar{\lambda}^{\infty,j_0,\bj}\geq0$, $\|\bar{\lambda}^{\infty,j_0,\bj}\|=1$, $\bar{w}^{\infty,j_0,\bj}\in \mathcal{N}_X(x^\infty)$ and $v_i^{\infty,j_0,\bj}\in\partial\zeta_i(x^\infty)$ for any $i\in\cI$. Dividing both sides of  (\ref{lg2222}) by $\|\lambda^{k,j_0,\bj}\|$ and taking limit as $\cK\ni k\to\infty$ yield
$\sum\limits_{i=1}^I\bar{\lambda}_i^{\infty,j_0,\bj}\left(\nabla \phi_i(x^\infty) +v^{\infty,j_0,\bj}_i-\nabla \psi_{i,j_i}(x^\infty)\right)+\bar{w}^{\infty,j_0,\bj}=0,$
which implies that
\begin{equation}\label{lg5555}
\sum\limits_{i=1}^I\bar{\lambda}_i^{\infty,j_0,\bj}\left(\nabla \phi_i(x^\infty)+v_i^{\infty,j_0,\bj}-\nabla \psi_{i,j_i}(x^\infty)\right)^Td=-(\bar{w}^{\infty,j_0,\bj})^Td\geq 0
\end{equation}
for any $d\in\mathcal{T}_X(x^\infty)$. On the other hand, since the PSCQ holds for $\Omega$ at $x^\infty$, it follows from  Definition \ref{slater} that there exists some $d_{\bj}\in \mathcal{T}_X(x^\infty)$ such that
\beq \label{pscq-i}
[\nabla \phi_i(x^\infty)+v_i^{\infty,j_0,\bj}-\nabla \psi_{i,j_i}(x^\infty)]^Td_{\bj}\leq \nabla \phi_i(x^\infty)^Td_{\bj}+\zeta_i'(x^\infty;d_{\bj})-\nabla \psi_{i,j_i}(x^\infty)^Td_{\bj}<0
\eeq
for all $i\in \cI_=(x^\infty)$. Notice that $\rho_k >0$, $\{\lambda^k/\rho_k\} \to 0$ and
\[
\{\hat{\phi}_i(x^{k,j_0,\bj};x^k)+\zeta_i(x^{k,j_0,\bj})-\ell_{\psi_{i,j_i}}(x^{k,j_0,\bj};x^k)\}_{k\in\cK} \to \phi_i(x^\infty)+\zeta_i(x^\infty)-\psi_i(x^\infty), \  \forall i\in \cI.
\]
 By these and (\ref{lam-def}), one can observe that $\bar{\lambda}_i^{\infty,j_0,\bj}=0$ for all $i\in \cI_<(x^\infty)$. Recall that $\bar{\lambda}^{\infty,j_0,\bj}\geq0$ and $\|\bar{\lambda}^{\infty,j_0,\bj}\|=1$. Hence, there exists some $\bar{i}\in\cI_=(x^\infty)$ such that $\bar{\lambda}_{\bar{i}}^{\infty,j_0,\bj}>0$. These together
with \eqref{pscq-i} imply that
$$\sum\limits_{i=1}^I\bar{\lambda}_i^{\infty,j_0,\bj}[\nabla \phi_i(x^\infty)+v_i^{\infty,j_0,\bj}-\nabla \psi_{i,j_i}(x^\infty)]^Td_{\bj}\le \bar{\lambda}_{\bar{i}}^{\infty,j_0,\bj}[\nabla \phi_{\bar{i}}(x^\infty)+v_{\bar{i}}^{\infty,j_0,\bj}-\nabla \psi_{\bar{i},j_{\bar{i}}}(x^\infty)]^Td_{\bj}<0,$$
which contradicts (\ref{lg5555}). Therefore, $\{\lambda^{k,j_0,\bj}\}_{k\in \cK}$ is bounded.

{We are now ready to complete the proof of statement (ii).} Indeed, by the boundedness of $\{\lambda^{k,j_0,\bj}\}_{k\in \cK}$, (\ref{lg2222}) and the fact that $\{v_i^{k,j_0,\bj}\}_{k\in \cK}$ is bounded, we immediately see that $\{w^{k,j_0,\bj}\}_{k\in \cK}$ is bounded. By the semicontinuity of $\partial\zeta_i(\cdot)$ and $\mathcal{N}_X(\cdot)$ (see \cite[Theorem 24.4]{Rockafellar}) and $\{x^{k,j_0,\bj}\}_{k\in \cK}\to x^\infty$, one can see that every accumulation point of $\{v^{k,j_0,\bj}_i\}_{k\in\cK}$ and $\{w^{k,j_0,\bj}\}_{k\in\cK}$ belongs to $\partial\zeta_i(x^\infty)$ and $\mathcal{N}_X(x^\infty)$, respectively. In view of these and (\ref{lg2222}), one can easily conclude that for each accumulation point $\lambda^{\infty,j_0,\bj}$ of $\{\lambda^{k,j_0,\bj}\}_{k\in \cK}$, there exist $v_i^{\infty,j_0,\bj}\in \partial\zeta_i(x^\infty)$ for every $i\in\{0,1,\ldots,I\}$ and $w^{\infty,j_0,\bj}\in \mathcal{N}_X(x^\infty)$ such that
\[\nabla \phi_0(x^\infty) +v^{\infty,j_0,\bj}_0-\nabla \psi_{0,j_0}(x^\infty)+\sum\limits_{i=1}^I\lambda_i^{\infty,j_0,\bj}
[\nabla \phi_i(x^\infty) +v^{\infty,j_0,\bj}_i-\nabla \psi_{i,j_i}(x^\infty)]+w^{\infty,j_0,\bj}=0.
\]
Moreover, using (\ref{lam-def})  and the facts that $\{\rho_k\}\to \infty$ and $\{\lambda^k/\rho_k\}\to0$, we can obtain that $\lambda^{\infty,j_0,\bj}\geq0$ and $\lambda^{\infty,j_0,\bj}_i=0$ for every $i\in\cI_<(x^\infty)$. Hence, statement (ii) holds.
\end{proof}

\begin{remark}
{ (i) Theorem \ref{alm-them-mtp} is established based on the assumption that $x^{k,j_0,\bj}\in X$ satisfies \eqref{dist-subdif} for any $j_0\in\cJ_0(x^\infty)$ and $\bj\in\cJ(x^\infty)$. Notice that $\cJ_0(x^\infty)\subseteq \cJ_{0,\epsilon}(x^k)$ and $\cJ(x^\infty)\subseteq \cJ_{\epsilon}(x^k)$, where $\{x^k\}_{k\in\cK}\to x^\infty$ for some subsequence $\cK$. Consequently, one sufficient condition for this assumption to hold is that for all sufficiently large $k$, $x^{k,j_0,\bj}$ satisfies \eqref{dist-subdif} for any $j_0\in\cJ_{0,\epsilon}(x^k)$ and $\bj\in\cJ_{\epsilon}(x^k)$.

(ii) From the proof of Theorem \ref{alm-them-mtp}, one can observe that the condition \eqref{dist-subdif} can be replaced by an alternative condition given in \eqref{app-opt}. Moreover,  $x^{k,j_0,\bj}\in X$ satisfying \eqref{app-opt} can be found by approximately solving $\min_{x\in X}\tQ_{\rho_k}(x^{k,j_0,\bj};x^k,\lambda^k,j_0,\mathbb{j})$ by mirror descent  or smoothing methods (e.g., see \cite{Nem83,Beck03,Nesterov05}).
}
\end{remark}

Before ending this section, we provide an example to illustrate the theoretical results of our AL method for solving problem \eqref{e1}.

\begin{example}\label{Example}
Consider the DC program
\beq\label{dc-example}
\begin{array}{ll}
\min\limits_{x\in \Re}& F(x)= |x|-\max\{6x,x\}\\
{\rm s.t.} &  2x-\max\{-x,x\}\le0.
\end{array}
\eeq

Clearly, it is a special case of \eqref{e1} with $\cI=\{1\}$, $\cJ_0=\cJ_1=\{1,2\}$, and
\[\ba{ll}
\phi_0(x)=0,\  \zeta_0(x)=|x|,\ \psi_{0,1}(x)=6x,\ \psi_{0,2}(x)=x,\\
\phi_1(x)=2x,\ \zeta_1(x)=0,\ \psi_{1,1}(x)=-x,\ \psi_{1,2}(x)=x.
\ea\]
We next apply the AL method, namely, Algorithm  \ref{aug-lag}  to solve problem \eqref{dc-example}. For convenience, we set $\epsilon=\infty$, $\alpha=1$, $\sigma=2$, $\lambda^0 = 0$, and
let $\rho_0>0$ be arbitrarily chosen for Algorithm  \ref{aug-lag}.  At the $k$th iteration,  we compute $x^k$ and update $\lambda^{k+1}$ and $\rho_{k+1}$ as follows.
\bi
\item[(i)]
We first compute
\beq\label{exam-sub}
x^{k,j_0,j_1}=\argmin\limits_{x\in\Re} {\bar F}_{\rho_k}(x;\lambda^k,j_0,j_1)
\eeq
for every $(j_0,j_1)\in\cJ_0\times\cJ_1$, where
$${\bar F}_{\rho_k}(x;\lambda^k,j_0,j_1)=\phi_0(x)+\zeta_0(x)-\psi_{0,j_0}(x)+\frac{1}{2\rho_k}[\lambda^k+\rho_k(\phi_1(x)+ \zeta_1(x)-\psi_{1,j_1}(x))]_+^2-\frac{(\lambda^k)^2}{2\rho_k}.$$
Then we set $x^k = x^{k,\hat{j}_0,\hat{j}_1}$ with $(\hat{j}_0,\hat{j}_1)$ given by
$$
(\hat{j}_0,\hat{j}_1)\in\Argmin\limits_{(j_0,j_1)} \{\tF_{\rho_k}(x^{k,j_0,j_1},\lambda^k)|(j_0,j_1)\in\cJ_0\times \cJ_1\},$$
where
$$ \tF_{\rho_k}(x,\lambda^k)=|x|-\max\{6x,x\}+\frac{1}{2\rho_k}[\lambda^k+\rho_k(2x-\max\{-x,x\})]_+^2-\frac{(\lambda^k)^2}{2\rho_k}.
$$
\item[(ii)] We update $\lambda^{k+1}$ and $\rho_{k+1}$ by
\[\lambda^{k+1}=[\lambda^k+\rho_k(2x^k-\max\{-x^k,x^k\})]_+,\quad \rho_{k+1}=\max\{2\rho_k,(\lambda^{k+1})^2\},\]
and let
\[
\lambda^{k,j_0,j_1}=[\lambda^k+\rho_k(\phi_1(x^{k,j_0,j_1})+ \zeta_1(x^{k,j_0,j_1})-\psi_{1,j_1}(x^{k,j_0,j_1}))]_+,\ \forall  (j_0,j_1)\in\cJ_0\times \cJ_1.
\]
\ei

By some simple calculations, one can find the expressions of $\{x^k\}$, $\{\lambda^k\}$, $\{x^{k,j_0,j_1}\}$ and $\{\lambda^{k,j_0,j_1}\}$, which are
presented in Table \ref{table}. Claim that $x^k$ and $x^{k,j_0,j_1}, \forall (j_0,j_1)\in\cJ_0\times\cJ_1$ satisfy \eqref{fun-value} and \eqref{dist-subdif} with
$\eta_k=0$, $\gamma_k=10/\rho_k$ and
\beq \label{exam-tQ}
\tQ_{\rho_k}(x;x^k,\lambda^k,j_0,j_1)={\bar F}_{\rho_k}(x;\lambda^k,j_0,j_1)+\frac{1}{2}(x-x^k)^2,
\eeq
which corresponds to \eqref{up-al-fun} with $L_0=1$ and $L_1=0$. Indeed, it is not hard to observe that
for every $(j_0,j_1)\in\cJ_0\times\cJ_1$ and $x\in\Re$,
\[
\tF_{\rho_k}(x^k,\lambda^k) \le \tF_{\rho_k}(x^{k,j_0,j_1},\lambda^k) \le  {\bar F}_{\rho_k}(x^{k,j_0,j_1};\lambda^k,j_0,j_1) \le {\bar F}_{\rho_k}(x;\lambda^k,j_0,j_1) \le \tQ_{\rho_k}(x;x^k,\lambda^k,j_0,j_1).
\]
It thus follows that $x^k$ satisfies \eqref{fun-value} with $\eta_k=0$.
In addition, one can see from Table \ref{table} that $|x^{k,j_0,j_1}|\leq 5/\rho_k$ and $|x^{k}|\leq 5/\rho_k$ for every $(j_0,j_1)\in\cJ_0\times\cJ_1$ and $k$. Also, from \eqref{exam-sub}, one has $0\in \partial {\bar F}_{\rho_k}(x^{k,j_0,j_1};\lambda^k,j_0,j_1)$. By these and \eqref{exam-tQ}, we obtain
$$ {\rm dist}\left(0, \partial [\tQ_{\rho_k}(x;x^{k},\lambda^k,j_0,j_1)]\big|_{x = x^{k,j_0,j_1}}\right)\le|x^{k}-x^{k,j_0,j_1}| \leq \frac{10}{\rho_k}, \ \forall (j_0,j_1)\in\cJ_0\times\cJ_1. $$
Hence, $x^{k,j_0,j_1}$ satisfies \eqref{dist-subdif} with $\gamma_k = 10/\rho_k$ for every $(j_0,j_1)\in\cJ_0\times\cJ_1$.

Notice that $\rho_k\rightarrow\infty$ as $k\rightarrow\infty$. Therefore, one can observe from Table \ref{table} that $\{x^k\}$ converges to $x^\infty=0$. It follows that $\cJ_0(x^\infty)=\cJ_1(x^\infty)=\{1,2\}$. Let $\lambda^{\infty,j_0,j_1}$ be any accumulation point of $\{\lambda^{k,j_0,j_1}\}$ for every $(j_0,j_1)\in\cJ_0(x^\infty)\times\cJ_1(x^\infty)$. By some simple calculations, one can verify that for every $(j_0,j_1)\in\cJ_0(x^\infty)\times\cJ_1(x^\infty)$,
\[\ba{l}
\lambda^{\infty,j_0,j_1}\geq0,\qquad \lambda^{\infty,j_0,j_1}\left[\phi_1(x^\infty)+\zeta_1(x^\infty)-\psi_{1,j_1}(x^\infty)\right]=0, \\
0  \in \nabla \phi_0(x^\infty)+\partial\zeta_0(x^\infty)-\nabla\psi_{0,j_0}(x^\infty)+\lambda^{\infty,j_0,j_1} [\nabla \phi_1(x^\infty)+\partial\zeta_1(x^\infty)-\nabla\psi_{1,j_1}(x^\infty)].
\ea\]
This result is indeed consistent with that in Theorem \ref{alm-them-mtp} since $x^\infty$ is a feasible point of \eqref{dc-example}, $\cI_=(x^\infty) = \{1\} \neq \emptyset$, and the PSCQ holds at $x^\infty$. The latter fact is due to $\nabla\phi_1(x^\infty)^Td+\zeta'_1(x^\infty;d)-\nabla\psi_{1,j_1}(x^\infty)^Td<0$ for every $j_1\in\cJ_1(x^\infty)$ and $d<0$.
\end{example}

\begin{table}[t]
\begin{center}
\caption{ Computational results of Algorithm \ref{aug-lag} for solving \eqref{dc-example}.}
\label{table}
\newcommand{\rb}[1]{\raisebox{1.5ex}[0pt]{#1}}
\begin{tabular}{c|c|c|c|c|c} \hline
Iteration & $(j_0,j_1)$ &$(1,1)$&$(1,2)$&$(2,1)$ &$(2,2)$\\
\hline\hline
& $x^{k,j_0,j_1}$&$5/(9\rho_k)$&$5/\rho_k$&$0$&$0$\\
\multirow{2}{*}{$k=0$}&$\lambda^{k,j_0,j_1}$&$5/3$&$5$&$0$&$0$\\
&$\tF_{\rho_k}(x^{k,j_0,j_1},\lambda^k)$&$-425/(162\rho_k)$&$-25/(2\rho_k)$&$0$&$0$\\
\cline{2-6}
&\multicolumn{5}{|c}{$x^k=5/\rho_k$,\quad $\lambda^{k+1}=5$}\\
\hline\hline
\multirow{4}{*}{$\ba{c}k=2m-1\\ (m=1,2,\ldots)\ea$}& $x^{k,j_0,j_1}$&$-8/(9\rho_k)$&$0$&$-13/(9\rho_k)$&$-3/\rho_k$\\
&$\lambda^{k,j_0,j_1}$&$7/3$&$5$&$2/3$&$2$\\
&$\tF_{\rho_k}(x^{k,j_0,j_1},\lambda^k)$&$-8/\rho_k$&$0$&$-169/(18\rho_k)$&$-13/(2\rho_k)$\\
\cline{2-6}
&\multicolumn{5}{|c}{$x^k=-13/(9\rho_k)$,\quad $\lambda^{k+1}=2/3$}\\
\hline\hline
\multirow{4}{*}{$\ba{c}k=2m\\ (m=1,2,\ldots)\ea$}& $x^{k,j_0,j_1}$&$1/(3\rho_k)$&$13/(3\rho_k)$&$0$&$0$\\
&$\lambda^{k,j_0,j_1}$&$5/3$&$5$&$2/3$&$2/3$\\
&$\tF_{\rho_k}(x^{k,j_0,j_1},\lambda^k)$&$-25/(18\rho_k)$&$-169/(18\rho_k)$&$0$&$0$\\
\cline{2-6}
&\multicolumn{5}{|c}{$x^k=13/(3\rho_k)$,\quad $\lambda^{k+1}=5$}\\
\hline
\end{tabular}
\end{center}
\end{table}

\section{Successive convex approximation method for penalty and AL subproblems}
\label{penalty-subprob}
{ An approximate solution of subproblems  \eqref{sbpp} and \eqref{al-pro1} satisfying \eqref{exawww10} and \eqref{fun-value}  is required in Algorithm \ref{PENM} and \ref{aug-lag}, respectively. Since these subproblems can be viewed as a special case of problem \eqref{exasubpm}, one can observe that to find these approximate solutions, it suffices to find an approximate solution $x_\eta$ of problem  (\ref{exasubpm}) satisfying that
\[
x_{\eta} \in X, \quad F_\rho(x_\eta) \le Q_\rho(x;x_\eta,j_0,\mathbb{j})+\eta, \ \ \forall x\in X, \ j_0\in \mathcal{J}_{0,\epsilon}(x_\eta), \ \mathbb{j}\in \mathcal{J}_{\epsilon}(x_\eta)
\]
for any given $\rho>0$, $\epsilon>0$ and $\eta>0$, where $Q_\rho$ is defined in \eqref{exaqappfun}. In what follows, we propose a successive convex approximation method to find such an approximate solution. The proposed method only solves a single convex problem in each iteration, while  the EDCA \cite{Pang} needs to solve a number of convex problems per iteration. It is therefore practically more efficient than the latter method.}



\begin{algorithm}
\label{I_SCA}
\normalfont
\mbox{}
\begin{itemize}
\item[0.] Input $x^0\in X$, $\eta>0$, $\epsilon>0$ and a sequence $\{\delta_t\}\subset\Re_+$ such that $\sum_{t=0}^\infty\delta_t^2 < \infty$. Set $\mathcal{B}^0\leftarrow\emptyset$ and $t\leftarrow 0$.
\item[1.]  Choose $(j_0,\mathbb{j})\in (\mathcal{J}_{0,\epsilon}(x^{t})\times \mathcal{J}_\epsilon(x^t))\setminus \mathcal{B}^t$, and find an approximate solution $x^{t,j_0,\mathbb{j}}$ of the problem
  \begin{equation}\label{exakkk111}
\min\limits_{x\in X} Q_\rho(x;x^t,j_0,\mathbb{j})
\end{equation}
 satisfying
\begin{equation}\label{termc}
x^{t,j_0,\mathbb{j}}\in X, \quad {\rm dist}\left(0, \partial [Q_\rho(x;x^t,j_0,\mathbb{j})+\iota_X(x)]\big|_{x=x^{t,j_0,\mathbb{j}}}\right) \le \delta_t,
\end{equation}
where $\iota_X$ is the indicator function of $X$.
\item[2.] If $F_\rho(x^t)-F_\rho(x^{t,j_0,\mathbb{j}})+\delta_t^2/(2L_0)>\eta$, set $x^{t+1}\leftarrow x^{t,j_0,\mathbb{j}}$, $\mathcal{B}^{t+1}\leftarrow\emptyset$, $t \leftarrow t+1$ and go to Step 1; otherwise, set $\mathcal{B}^t\leftarrow\mathcal{B}^t\cup\{(j_0,\mathbb{j})\}$ and go to Step 3.
\item[3.] If $\mathcal{J}_{0,\epsilon}(x^{t})\times \mathcal{J}_\epsilon(x^t)=\mathcal{B}^t$, stop; otherwise, go to Step 1.
\end{itemize}
{\bf End.}
\end{algorithm}

\begin{remark}
(i) In contrast with \eqref{exasubpm}, problem \eqref{exakkk111} has a simpler objective function and it can be efficiently solved for many $X$ and $\zeta_i$'s. For example, when $X$ is a polyhedral set or more generally a conic quadratic representable set, and $\zeta_i$'s are polyhedral functions or more generally conic quadratic representable functions, problem \eqref{exakkk111} can be reformulated as a conic quadratic program, which can be efficiently solved by interior point methods.

(ii) As seen from the proof of Theorem \ref{SCA_C} below, the condition \eqref{termc} can be replaced by {an alternative condition:}
\[
x^{t,j_0,\mathbb{j}}\in X, \quad Q_\rho(x^{t,j_0,\mathbb{j}};x^t,j_0,\mathbb{j})\leq  Q_\rho(x;x^t,j_0,\mathbb{j})+ \delta_t^2/(2L_0), \quad \forall x\in X.
\]
That is, $x^{t,j_0,\mathbb{j}}$ is a $\delta_t^2/(2L_0)$-optimal solution
of problem \eqref{exakkk111}, which can be found by mirror descent  or smoothing methods (e.g., see \cite{Nem83,Beck03,Nesterov05}).
\end{remark}

We now establish some convergence results for Algorithm \ref{I_SCA}. 


\begin{theorem}\label{SCA_C}
Assume that the function $F_\rho$ is bounded below on $X$.\footnote{It can be seen that this assumption holds if $F$ is bounded below on $X$.} Then Algorithm \ref{I_SCA} terminates in finitely many iterations, that is, there exists an integer $\hat{t} \ge 0$ such that
\beq\label{K_Iq}
F_\rho(x^{\hat{t}})-F_\rho(x^{\hat{t},j_0,\mathbb{j}})+\delta_{\hat{t}}^2/(2L_0)\leq\eta, \quad
\forall j_0\in \mathcal{J}_{0,\epsilon}(x^{\hat{t}}), \forall \mathbb{j}\in \mathcal{J}_\epsilon(x^{\hat{t}}),
\eeq
where $x^{\hat{t}}$ and $x^{\hat{t},j_0,\mathbb{j}}$ are generated by Algorithm \ref{I_SCA}
for all $j_0\in \mathcal{J}_{0,\epsilon}(x^{\hat{t}})$ and $\mathbb{j}\in \mathcal{J}_\epsilon(x^{\hat{t}})$.
Moreover, for any $j_0\in\mathcal{J}_{0,\epsilon}(x^{\hat{t}})$ and $\mathbb{j}\in\mathcal{J}_{\epsilon}(x^{\hat{t}})$, it holds
\beq\label{xk_est}
F_\rho(x^{\hat{t}})\leq Q_\rho(x;x^{\hat{t}},j_0,\mathbb{j})+\eta,\quad\forall x\in X.
\eeq
\end{theorem}

\begin{proof}
Suppose for contradiction that Algorithm \ref{I_SCA} does not terminate in finitely many iterations. Let $\{x^t\}$ be the sequence generated by Algorithm \ref{I_SCA}. Then it follows from Steps 2 and 3 of Algorithm \ref{I_SCA} that
\[
F_\rho(x^t)-F_\rho(x^{t+1})+\delta_t^2/(2L_0)>\eta,\quad\forall t \ge 0,
\]
which implies that
\[
F_\rho(x^t)<F_\rho(x^0)+\frac{1}{2L_0}\sum\limits_{i=0}^{t-1}\delta_i^2-t\eta,\quad \forall t \ge 1.
\]
By this and $\sum_{t=0}^\infty\delta_t^2 < \infty$, one can obtain that $\lim\limits_{t\to\infty}F_\rho(x^t)=-\infty$, which, together with $\{x^t\} \subset X$, contradicts the assumption that $F_\rho$ is bounded below on $X$. Hence, Algorithm \ref{I_SCA} terminates after finitely many iterations, which implies that $\mathcal{B}^{\hat{t}}=\mathcal{J}_{0,\epsilon}(x^{\hat{t}})\times \mathcal{J}_\epsilon(x^{\hat{t}})$ and \eqref{K_Iq} hold for some $\hat t \ge 0$.

Given that \eqref{K_Iq} holds for some $\hat{t}$, we next show that \eqref{xk_est} holds for any $j_0\in \mathcal{J}_{0,\epsilon}(x^{\hat{t}})$ and $\mathbb{j}\in\mathcal{J}_\epsilon(x^{\hat{t}})$. To this end, let $j_0\in \mathcal{J}_{0,\epsilon}(x^{\hat{t}})$ and $\mathbb{j}\in\mathcal{J}_\epsilon(x^{\hat{t}})$ be arbitrarily chosen. One can observe from (\ref{hphi}) and (\ref{exaqappfun}) that $Q_{\rho}(\cdot;x^{\hat{t}},j_0,\mathbb{j})$ is strongly convex on $X$ with modulus $L_0>0$. Since $x^{{\hat{t}},j_0,\bj}$ satisfies (\ref{termc}),  by a similar argument as for deriving \eqref{app-opt}, we have
\[Q_{\rho}(x^{{\hat{t}},j_0,\bj};x^{\hat{t}},j_0,\mathbb{j})\leq Q_{\rho}(x;x^{\hat{t}},j_0,\mathbb{j})+\delta_{\hat{t}}^2/(2L_0),\quad \forall x\in X,\]
which yields
\[
F_\rho(x^{{\hat{t}},j_0,\mathbb{j}})\leq   Q_\rho(x^{{\hat{t}},j_0,\mathbb{j}};x^{\hat{t}},j_0,\mathbb{j})\leq  Q_\rho(x;x^{\hat{t}},j_0,\mathbb{j})+ \delta_{\hat{t}}^2/(2L_0),\quad \forall x\in X.
\]
This together with \eqref{K_Iq} implies that \eqref{xk_est} holds as desired.
\end{proof}

\setcounter{equation}{0}
\section{Numerical results}
\label{Num-Res}

In this section we conduct some numerical experiments to test the performance of our proposed methods, namely, the penalty method (PM) in Algorithm \ref{PENM} and the augmented Lagrangian method (ALM) in Algorithm \ref{aug-lag}, and compare them with two closely related methods proposed in \cite[Section 6]{Pang}, which are an enhanced DCA (EDCA) and an exact penalty method (EPM). {For convenience, we use PM1 and PM2 to stand for the PM with $p=1$ and $2$, respectively. The subproblems \eqref{sbpp} and \eqref{al-pro1} of PM and ALM are solved by Algorithm \ref{I_SCA}.
Also, the EDCA requires a feasible initial point to start while the other methods do not.
We will compare these methods numerically below. All the methods are coded in Matlab and all the computations are performed on a Dell laptop with an Intel Core i7-1065G7 CPU and 16 GB of RAM.} 

In the first experiment, we apply the aforementioned methods to the following optimization problem with two structured DC constraints:
\begin{equation}
\label{eq:exp-1}
\min_{x\in\Re^n} \left\{ \phi_0(x) \ | \ \phi_i(x) - \psi_i(x) \leq 0, \ \forall i=1,2\right\}, \quad \mbox{where} \ \ \psi_i(x) = \max\{ \psi_{i,1}(x), \psi_{i,2}(x)\}, \ i=1,2.
\end{equation}
The functions $\phi_i$'s and $\psi_{i,j}$'s in \eqref{eq:exp-1} are convex quadratic functions, namely,
\begin{align*}
\phi_0(x) & = x^TQx + q^Tx, \\
\phi_i(x) & = x^TA_ix + a_i^Tx + c_i \quad i=1,2, \\
\psi_{i,j}(x) & = x^TB_{i,j}x + b_{i,j}^Tx + d_{i,j} \quad i=1,2, \ j=1,2,
\end{align*}
where $Q,A_i,B_{i,j}\in\Re^{n\times n}$ are positive semidefinite matrices, $q,a_i,b_{i,j}\in\Re^n$, and $c_i,d_{i,j}\in\Re$ for all $i$ and $j$. It is clear that \eqref{eq:exp-1} is a special case of problem \eqref{e1}.

In this experiment, we set $\epsilon=0.01$ in all the above methods, $\rho_0 = 0.1$ and $\sigma=2$ in PM, ALM and EPM, $\alpha=1.05$ in ALM, and $\eta_k=10^{-k-3}$ for all $k$ in PM and ALM. In Algorithm \ref{I_SCA}, we set $\delta_t=10^{-t-1}$ for all $t$. For this test example,
 when applied to solve the subproblems of PM2 and ALM, the subproblem \eqref{exakkk111} of Algorithm \ref{I_SCA} is smooth and solved by a nonmonotone gradient method \cite{Ra97}. On the other hand, when applied to solve the subproblems of PM1,
the subproblem \eqref{exakkk111} of Algorithm \ref{I_SCA} is nonsmooth and solved by CVX\footnote{CVX is a Matlab package for solving convex programs; see \url{cvxr.com/cvx/}.}.
In addition, the penalty subproblem of EPM is solved as described in \cite{Pang}. In particular, when applied to \eqref{eq:exp-1}, the penalty subproblem of EPM is in the form of
\begin{equation}
\label{eq:epm-sub-o}
\min_{x\in\Re^n} \phi_0(x) + \rho\left[\max\{\phi_1(x) - \psi_1(x), \phi_2(x) - \psi_2(x)\}\right]_+,
\end{equation}
where $\rho>0$ is the penalty parameter. As described in \cite{Pang}, we first rewrite \eqref{eq:epm-sub-o} as
\begin{equation}
\label{eq:epm-sub}
\min_{x\in\Re^n} \underbrace{\phi_0(x) + \rho\max\{\psi_1(x) + \psi_2(x), \phi_1(x) + \psi_2(x), \phi_2(x) + \psi_1(x)\}}_{\mathrm{convex}} - \underbrace{\rho(\psi_1(x) + \psi_2(x))}_{\mathrm{convex}}
\end{equation}
and then apply \cite[Algorithm 1]{Pang} to solve \eqref{eq:epm-sub}. Moreover, the convex problems arising in each iteration of \cite[Algorithm 1]{Pang} are solved by CVX. While in theory the EPM requires an exact D-stationary point of its penalty subproblem \eqref{eq:epm-sub}, in our experiment we terminate \cite[Algorithm 1]{Pang} when the norm of the difference of two consecutive iterates generated by it is less than $10^{-6}$.
Finally, the subproblems of EDCA are constrained convex programs and solved by CVX.

We randomly generate $4$ instances for problem \eqref{eq:exp-1} with $n=50,100,250, 500$, respectively, each of which is generated as follows. Given a positive integer $n$, we first generate a vector $d\in\Re^n$, whose entries are randomly chosen from a uniform distribution on $[0,20]$. We then generate a matrix $\tilde{U}\in\Re^{n\times n}$ with entries randomly chosen from the standard normal distribution, compute an orthogonal basis $U$ for the range space of $\tilde{U}$, and set $Q = U \mathrm{Diag}(d) U^T$. The matrices $A_i$ and $B_{i,j}$ for all $i$ and $j$ are randomly generated in the same manner as $Q$. In addition, we generate the vectors $q, a_i,b_{i,j}\in\Re^n$ for all $i$ and $j$ with entries randomly chosen from the standard normal distribution.
Also, we randomly choose $c_i$ and $d_{i,j}$ for all $i$ and $j$ from the standard normal distribution.

For each instance, we perform $10$ runs of all the tested methods as described above. In each run, we first randomly generate a point $x^0$, whose entries are randomly chosen from the standard normal distribution. We then run PM, ALM, and EPM with the same initial point $x^0$. Moreover, if $x^0$ is feasible for \eqref{eq:exp-1}, we run the EDCA with the initial point $x^0$; if not, we repeat generating $\tilde{x}^0$ with entries randomly chosen from the standard normal distribution until a feasible point $\tilde{x}^0$ is found and then run the EDCA with the initial point $\tilde{x}^0$. We terminate all the tested methods once 
$\|x^{k+1}-x^{k}\|/\|x^{k+1}\|\leq10^{-5}$ holds for some $k$, where $x^k$ and $x^{k+1}$ are the approximate solutions obtained at the $k$th and $(k+1)$th iterations of each method, respectively.

The computational results averaged over each group of 10 runs with same $n$ are presented in Table \ref{tab:0}, which consists of four subtables. In detail, the parameter $n$ is listed in the first column. For each $n$, the objective value at the solutions produced by all the tested methods, averaged over 10 runs, is given in Table \ref{tab:0}(a), and the outer iteration number, the total number of convex subproblems solved, and the CPU time (in seconds) averaged over 10 runs are given in Tables \ref{tab:0}(b), \ref{tab:0}(c), and \ref{tab:0}(d), respectively. From Table \ref{tab:0}(a), one can see that the solutions produced by all the above methods have about the same objective value. Also, as seen from Table \ref{tab:0}(b), EPM takes less outer iterations than PM1, PM2 and ALM, while EDCA takes much more outer iterations. In addition, one can see from Table \ref{tab:0}(c) that PM2 generally solves more convex subproblems than
the other methods, while PM1 solves less convex subproblems. Besides,
Table \ref{tab:0}(d) shows that PM1, PM2 and ALM are much faster than EPM and EDCA, which is mainly because the convex subproblems arising in the latter two methods are more sophisticated and solved by CVX. Also, PM2 and ALM are much faster than PM1, which is largely due to the fact that the convex subproblems arising in PM2 and ALM are smooth and solved by a gradient method, while the ones arising in PM1 are nonsmooth and solved by CVX.

\begin{table}[t]
\centering
	\caption{Computational results for solving problem \eqref{eq:exp-1}}
	\label{tab:0}
	\begin{subtable}[t]{\linewidth}
	\centering
		\begin{tabular}{|c ||c c c c c|}
			\hline
			$n$ & PM1   & PM2  & ALM & EPM & EDCA \\
			\hline	
			50 & -1.901 & -1.901 & -1.901 & -1.901 & -1.901\\
			100 & -3.825 & -3.825 & -3.825 & -3.825 & -3.825 \\
			250 & -9.061 & -9.061 & -9.061 & -9.060 & -9.061 \\
			500 & -22.390 & -22.390 & -22.390 & -22.390 & -22.390\\
			\hline	
	\end{tabular}
	\caption{Results for objective value}
	\end{subtable}

\vspace{0.2in}	
	
	\begin{subtable}[t]{\linewidth}
	\centering
		\begin{tabular}{|c ||c c c c c|}
			\hline
			$n$ & PM1   & PM2  & ALM & EPM & EDCA \\
			\hline
			50 & 3.0 & 12.0& 4.1& 2.0& 20.6\\
			100 & 3.0 & 9.2& 3.4& 2.6& 17.4\\
			250 & 3.1 & 10.0& 4.0 &2.0 &24.0 \\
			500 & 2.5& 2.0& 2.0& 2.0& 12.8\\
			\hline	
	\end{tabular}
	\caption{Results for outer iteration number}
	\end{subtable}

\vspace{0.2in}	
		
	\begin{subtable}[t]{\linewidth}
	\centering
		\begin{tabular}{|c ||c c c c c|}
			\hline
			$n$ & PM1   & PM2  & ALM & EPM & EDCA \\
			\hline
			50 & 13.4 & 84.7& 20.9& 32.9& 21.0\\
			100 & 13.3 & 50.0& 18.9& 34.4& 18.2\\
			250 & 15.2 & 71.2& 23.8 &28.0 &24.8 \\
			500 & 10.7& 11.7& 11.7& 33.6& 13.2\\
			\hline	
	\end{tabular}
	\caption{Results for number of convex subproblems solved}
	\end{subtable}

\vspace{0.2in}	

	\begin{subtable}[t]{\linewidth}
	\centering
		\begin{tabular}{|c ||c c c c c|}
			\hline
			$n$ & PM1   & PM2 & ALM & EPM & EDCA \\
			\hline
			50 & 4.4 & 0.4 & 0.1 & 21.5 & 10.1 \\
			100 & 4.6 & 0.2 & 0.1 & 28.3 & 10.0 \\
			250 & 12.8 & 4.6 & 1.1 & 149.4 & 50.4 \\
			500 & 33.3 & 2.9 & 2.9 & 707.3 & 170.4 \\
			\hline	
		\end{tabular}
		\caption{Results for CPU time}
	\end{subtable}
\end{table}

In the second experiment, we apply PM, ALM, EPM and EDCA to the following DC program:
\begin{equation}
\label{eq:exp-2}
\min_{x\in\Re^n} \left\{ \|Ax - b\|^2 \ | \ \|x\|_1 - h(x) \leq sK\right\},
\end{equation}
where $A\in\Re^{m\times n}$, $b\in\Re^m$, $1\leq K\leq n$ is an integer, $s>0$, and $h(x)$ is defined as
\[h(x)=\sum\limits_{i=1}^n\max\{x_i-s,0,-x_i-s\}.\]
Problem \eqref{eq:exp-2} arises in applications such as sparse signal recovery (e.g., see~\cite{Zheng}). It is clear that \eqref{eq:exp-2} is a special case of problem \eqref{e1}.

In this experiment, the same parameters as in the first experiment are chosen for all the tested methods, except $\epsilon = 0.01, 0.05$. For the test problem \eqref{eq:exp-2}, when applied to solve the subproblems of PM2 and ALM, the subproblem \eqref{exakkk111} of Algorithm \ref{I_SCA} is solved by a first-order method. On the other hand, when applied to solve the subproblems of PM1, the subproblem \eqref{exakkk111} of Algorithm \ref{I_SCA} is nonsmooth and solved by CVX.
The penalty subproblems of EPM are solved by \cite[Algorithm 1]{Pang} and the convex problems arising in each iteration of \cite[Algorithm 1]{Pang} are solved by CVX. We terminate \cite[Algorithm 1]{Pang} when the norm of the difference of two consecutive iterates generated by it is less than $10^{-6}$. Also, the subproblems of the EDCA are solved by CVX. We terminate all the tested methods once
$\|x^{k+1}-x^{k}\|/\|x^{k+1}\|\leq10^{-5}$, where $x^k$ and $x^{k+1}$ are the approximate solutions obtained at the $k$th and $(k+1)$th iterations of each method, respectively.

We choose $(m,n) = (2^8, 2^{10})$, $K = 20, 30, 40$, and set $s = 0.1$
in our experiment. For each $K$, we randomly generate $10$ instances of problem \eqref{eq:exp-2} in a similar manner as described in~\cite{LL18}. Given $K$, we first randomly generate a $K$-sparse vector $x^*\in\Re^n$. Specifically, we randomly choose $K$ numbers from $\{1,2,\dots,n\}$ as the support for $x^*$ and randomly choose the nonzero entries of $x^*$ from $\{-1,1\}$ with equal probability. We then generate the $m\times n$ data matrix $A$, whose entries are randomly chosen from the standard normal distribution. Finally we orthonormalize the rows of $A$ and set $b = Ax^* + \xi$, where the entries of $\xi\in\Re^n$ are drawn from a normal distribution with mean $0$ and variance $10^{-3}$. For each such instance, we apply the above methods to solve \eqref{eq:exp-2}. Since EDCA needs a feasible point of \eqref{eq:exp-2} to start, we first solve the following convex program by CVX:
$$ \min_{x\in\Re^n} \left\{ \|Ax - b\|^2 \ | \ \|x\|_1 \leq sK \right\}, $$
whose optimal solution $\tilde{x}^0$ must be feasible for \eqref{eq:exp-2}. We then run PM, ALM, EPM and  EDCA with $\tilde{x}^0$ as the initial point.
For all the above methods, we compute the relative error of the final iterate $\tilde{x}$ produced by them according to $\mathrm{rel\_err} = \|\tilde{x} - x^*\|/\|x^*\|$, which evaluates how well the sparse vector $x^*$ is recovered by $\tilde{x}$.

The computational results of this experiment are presented in Tables \ref{tab:1} and \ref{tab:2}. In detail, the parameter $K$ is listed in the first column. For each $K$, the objective value and the relative error at the solutions produced by all the tested methods, averaged over $10$ runs, are given in Tables \ref{tab:1}(a) and \ref{tab:1}(b), respectively, and the outer iteration number, the number of convex subproblems solved and the CPU time (in seconds) averaged over $10$ runs are given in Tables \ref{tab:2}(a), \ref{tab:2}(b) and \ref{tab:2}(c), respectively. From Table \ref{tab:1}, one can see that the objective value and the relative error at the solutions produced by the tested methods are about the same except the cases $K=30$ and $40$, for which those given by EDCA are larger. 
Also, from Table \ref{tab:2}, we observe that the outer iteration number, the number of convex subproblems solved, and the CPU time taken by PM1, PM2, and ALM are almost same for different $\epsilon$. However, as $\epsilon$ increases, EPM and EDCA solve many more convex subproblems and thus take much more CPU time. 
In addition,  ALM and PM2 are much faster than PM1, EPM and EDCA because the convex subproblems in ALM and PM2 have a simpler structure than those in the other methods and are solved more cheaply.


\begin{table}[htbp]
\centering
	\caption{Computational results for solving problem \eqref{eq:exp-2} }
	\label{tab:1}
	\begin{subtable}[t]{\linewidth}
	\centering
		\begin{tabular}{|c ||c c c c c||c c c c c|}
			\hline
             & \multicolumn{5}{|c||}{$\epsilon=0.01$} & \multicolumn{5}{|c|}{$\epsilon=0.05$} \\
             \cline{2-11}
			$K$ & PM1   & PM2  & ALM &EPM & EDCA& PM1   & PM2  & ALM & EPM & EDCA\\
			\hline
			20 &2.3e-04 & 2.3e-04 & 2.3e-04&  2.3e-04  & 2.3e-04&2.3e-04 & 2.3e-04 & 2.3e-04&  2.3e-04  & 2.3e-04\\
			30 &2.3e-04 &2.3e-04 & 2.3e-04 &2.3e-04 & 6.9e-02&2.3e-04  &2.3e-04 & 2.3e-04&2.3e-04 & 3.9e-02\\
			40 &2.1e-04 &2.1e-04 & 2.1e-04 &2.1e-04 & 1.8e-01&2.1e-04 &2.1e-04 & 2.1e-04 &2.1e-04 & 9.6e-02\\
			\hline	
	\end{tabular}
	\caption{Results for objective value}
	\end{subtable}
			
	\begin{subtable}[t]{\linewidth}
	\centering
			\begin{tabular}{|c ||c c c c c||c c c c c| }
			\hline
             & \multicolumn{5}{|c||}{$\epsilon=0.01$} & \multicolumn{5}{|c|}{$\epsilon=0.05$} \\
             \cline{2-11}
			$K$ & PM1   & PM2  & ALM &EPM  & EDCA& PM1   & PM2  & ALM &  EPM & EDCA\\
			\hline
			20 &2.0e-03 &2.0e-03 & 2.0e-03 &2.0e-03  & 2.0e-03&2.0e-03 &2.0e-03 & 2.0e-03 &2.0e-03  & 2.0e-03\\
			30  &1.9e-03 &1.9e-03 & 1.9e-03 &1.9e-03 & 5.3e-02&1.9e-03 &1.9e-03 & 1.9e-03 &1.9e-03 & 2.9e-02\\
			40  &2.1e-03 &2.1e-03 & 2.1e-03 &2.1e-03 & 1.1e-01&2.1e-03 &2.1e-03 & 2.1e-03 &2.1e-03 & 7.3e-02\\
			\hline	
		\end{tabular}
		\caption{Results for relative error}
	\end{subtable}
\end{table}

\begin{table}[htbp]
\centering
	\caption{Computational results for solving problem \eqref{eq:exp-2} }
	\label{tab:2}

	\begin{subtable}[t]{\linewidth}
	\centering
			\begin{tabular}{|c ||c c c c c||c c c c c|}
			\hline
             & \multicolumn{5}{|c||}{$\epsilon=0.01$} & \multicolumn{5}{|c|}{$\epsilon=0.05$} \\
             \cline{2-11}
			$K$ & PM1   & PM2  & ALM &EPM & EDCA & PM1   & PM2  & ALM &EPM & EDCA \\
			\hline
			20&10.0 &10.0 & 9.4 &2.0 & 6.1 &10.0 &10.0 & 9.4 &2.0 & 4.9 \\
			30 &9.8 &9.9 & 9.0 &2.0 & 7.4 &9.8 &9.9 & 9.0 &2.0  & 5.6 \\
			40 &9.6 &9.5 & 9.3 &2.2 & 8.0 &9.5 &9.5 & 9.3 &2.2  & 7.0 \\
			\hline	
	\end{tabular}
	\caption{Results for outer iteration number}
	\end{subtable}

	\begin{subtable}[t]{\linewidth}
	\centering
			\begin{tabular}{|c ||c c c c c||c c c c c|}
			\hline
             & \multicolumn{5}{|c||}{$\epsilon=0.01$} & \multicolumn{5}{|c|}{$\epsilon=0.05$} \\
             \cline{2-11}
			$K$ & PM1   & PM2  & ALM &EPM & EDCA & PM1   & PM2  & ALM &EPM & EDCA \\
			\hline
			20&88.9 &116.8 & 108.8 &5.9 & 10.2 &88.9 &116.8 & 108.8 &46.5 & 66.2 \\
			30 &103.1 &132.0 & 117.0 &6.1 & 11.2 &101.4 &132.0 & 117.6 &21.7 & 88.0\\
			40 &134.9 &152.5 & 144.4 &8.9 & 17.0 &136.7 &152.5 & 144.4 &299.1  & 449.8 \\
			\hline	
	\end{tabular}
	\caption{Results for number of convex subproblems solved}
	\end{subtable}

	\begin{subtable}[t]{\linewidth}
	\centering
		\centering
			\begin{tabular}{|c ||c c c c c||c c c c c|}
			\hline
             & \multicolumn{5}{|c||}{$\epsilon=0.01$} & \multicolumn{5}{|c|}{$\epsilon=0.05$}\\
             \cline{2-11}
			$K$ & PM1   & PM2   &  ALM  &EPM & EDCA & PM1   & PM2   &  ALM  &EPM & EDCA \\
			\hline
			20&28.6 & 0.07 & 0.03 &36.7 & 32.3 &29.3 & 0.06 & 0.03 &414.6 & 203.6 \\
			30 &32.4 & 0.07 & 0.03 &41.2 & 35.1 &31.8 & 0.06 & 0.03 &201.8 & 276.7 \\
			40 &49.3 &0.10 & 0.04 &78.3 & 69.9 &50.1 &0.08 & 0.04 &2872.6  & 1482.4 \\
			\hline	
		\end{tabular}
		\caption{Results for CPU time}
	\end{subtable}
\end{table}

\section{Concluding remarks}
\label{conclude}
The current development of this paper is based on the assumption that the second convex component of the objective and constraints is the supremum of finitely many convex smooth functions. It is worthy of a further research whether it can be extended to the case where the second convex component is the supremum of infinitely many convex smooth functions. 

\appendix
\section{Proof of Theorem \ref{thm:KKT}} \label{proof-KKT}

In this section we provide a proof of Theorem \ref{thm:KKT}. Before proceeding, we establish a technical lemma as follows.

\begin{lemma}
\label{polar-cone}
Let $\bar{x}\in\Omega$ be such that $\mathcal{I}_=(\bar{x})\neq\emptyset$, and let
\beq\label{PJ-cone}
P_\bj(\bx)=\left\{\sum\limits_{i\in\cI_=(\bx)}\lambda_i[\nabla \phi_i(\bx)+v_i-\nabla\psi_{i,j_i}(\bx)]+w:\ \lambda_i\geq0,  v_i\in\partial\zeta_i(\bx),  w\in\mathcal{N}_X(\bx)\right\}
\eeq
for any $\bj=(j_1,\ldots,j_I)\in\cJ(\bx)$. Then $[P_\bj(\bx)]^\o=C_\bj(\bx)$ for any $\bj\in\cJ(\bx)$, where $C_\bj(\bx)$ is defined in \eqref{lin-cone} and ${\cal S}^\o$ denotes the polar cone of any cone $\cal S$. Assume further that the  {\rm PSCQ} holds for $\Omega$ at $\bar{x}$, or that $X$ is a polyhedral set, $\phi_{i}$ is affine and $\zeta_i$ is piecewise affine on $X$ for every $i\in\mathcal{I}_=(\bar{x})$. Then $P_\bj(\bx)$ is a nonempty closed convex cone and $P_\bj(\bx)=[C_\bj(\bx)]^\o$ for any $\bj\in\cJ(\bx)$.
\end{lemma}

\begin{proof}
Let $\bj=(j_1,\ldots,j_I)\in\cJ(\bx)$ be arbitrarily chosen. We first prove $C_\bj(\bx)\subseteq [P_\bj(\bx)]^\o$. To this end, let $d\in C_\bj(\bx)$ be arbitrarily chosen. It then follows from \eqref{lin-cone} that $d\in \cT_X(\bx)$ and $\nabla \phi_i(\bx)^Td+\zeta'(\bx;d)-\nabla\psi_{i,j_i}(\bx)^Td\leq 0$ for each $i\in\cI_=(\bx)$.  By $d\in \cT_X(\bx)$, one has that $w^Td\leq0$ for any $w\in\mathcal{N}_X(\bx)$. Clearly, one also has $v_i^Td\le\zeta'(\bx;d)$ for any $v_i\in\partial\zeta_i(\bx)$. In view of these, we have that for any $\lambda_i\ge0$ with $i\in\cI_=(\bx)$,
\[d^T\left(\sum\limits_{i\in\cI_=(\bx)}\lambda_i[\nabla \phi_i(\bx)+v_i-\nabla\psi_{i,j_i}(\bx)]+w\right)\le\sum\limits_{i\in\cI_=(\bx)}\lambda_i[\nabla \phi_i(\bx)^Td+\zeta'(\bx;d)-\nabla\psi_{i,j_i}(\bx)^Td]\leq0.
\]
Hence, $d\in[P_\bj(\bx)]^\o$, which leads to $C_\bj(\bx)\subseteq [P_\bj(\bx)]^\o$. We next prove that $[P_\bj(\bx)]^\o \subseteq C_\bj(\bx)$. Let $d\in[P_\bj(\bx)]^\o$ be arbitrarily chosen. Notice from \eqref{PJ-cone} that $\mathcal{N}_X(\bx)\subseteq P_\bj(\bx)$. It thus follows that $w^Td\leq0$ for all $w\in \mathcal{N}_X(\bx)$, which implies $d\in \mathcal{T}_X(\bx)$. In addition, observe that $\nabla \phi_i(\bx)+v_i-\nabla \psi_{i,j_i}(\bar{x})\in P_\bj(\bx)$ for every $i\in\cI_=(\bx)$ and $v_i\in\partial\zeta_i(\bx)$, which along with $d\in[P_\bj(\bx)]^\o$ implies $[\nabla \phi_i(\bx)+v_i-\nabla \psi_{i,j_i}(\bar{x})]^Td\le0$. It follows
\[\nabla\phi_i(\bx)^Td+\zeta_i'(\bx;d)-\nabla \psi_{i,j_i}(\bx)^Td=\max_{v_i\in\partial\zeta_i(\bx)}\{[\nabla \phi_i(\bx)+v_i-\nabla \psi_{i,j_i}(\bar{x})]^Td\}\leq0,\ \forall i\in\cI_=(\bx).\]
Hence, $d\in C_\bj(\bx)$, which yields $[P_\bj(\bx)]^\o\subseteq C_\bj(\bx)$. This together with
$C_\bj(\bx)\subseteq [P_\bj(\bx)]^\o$ implies $[P_\bj(\bx)]^\o=C_\bj(\bx)$.

We next  prove that $P_\bj(\bx)$ is a closed convex cone under the assumption that
the  {\rm PSCQ} holds for $\Omega$ at $\bar{x}$, or that $X$ is a polyhedral set, $\phi_{i}$ is affine and $\zeta_i$ is piecewise affine on $X$ for every $i\in\mathcal{I}_=(\bar{x})$.

Firstly,  we assume that the  {\rm PSCQ} holds for $\Omega$ at $\bar{x}$. Clearly, $P_\bj(\bx)$ is a convex cone. Suppose for contradiction that $P_\bj(\bx)$ is not closed. Then there are some $u\not\in P_\bj(\bx)$ and some sequences $\{\lambda_i^k\}$, $\{v^k_i\}\subseteq\partial\zeta_i(\bx)$ and $\{w^k\}\subseteq\mathcal{N}_X(\bx)$ with $\lambda_i^k\geq0$
such that
\beq \label{u-eqn}
\lim\limits_{k\to\infty}\sum\limits_{i\in\cI_=(\bx)}\lambda_i^k[\nabla \phi_i(\bx)+v^k_i-\nabla\psi_{i,j_i}(\bx)]+w^k=u.
\eeq
Notice that $\partial\zeta_i(\bx)$ is compact for all $i\in\cI_=(\bx)$. It follows that $\{v^k_i\}$ is bounded. In view of these and the closedness of $\mathcal{N}_X(\bx)$, it is not hard to observe that there exists some $\hat{i}\in\cI_=(\bx)$ such that $\{\lambda^k_{\hat{i}}\}$ is unbounded (otherwise, one would have from \eqref{u-eqn}
 that $\{w^k\}$ is bounded, which along with \eqref{u-eqn}, the boundedness of $\{v^k_i\}$ and the closedness of $\partial\zeta_i(\bx)$ and $\mathcal{N}_X(\bx)$ implies that $u\in P_\bj(\bx)$). By passing to a subsequence if necessary, we can assume that $\{\lambda^k_{\hat{i}}\} \to \infty$.  Since the  {\rm PSCQ} holds for $\Omega$ at $\bar{x}$, there exists some vector $d_\bj\in\cT_X(\bx)$ such that $\nabla \phi_i(\bx)^Td_\bj+\zeta'_i(\bar{x};d_\bj)-\nabla \psi_{i,j_i}(\bar{x})^Td_\bj<0$ for all $i\in\cI_=(\bx)$. Also, by $\{v^k_i\}\subseteq\partial\zeta_i(\bx)$ and $\{w^k\}\subseteq\mathcal{N}_X(\bx)$, one has that  $d_\bj^T v^k_{i}\le\zeta'_{i}(\bar{x};d_\bj)$ and $d_\bj^T w^k\le 0$ for all $k$. Using these, \eqref{u-eqn} , $\lambda^k_i \ge 0$ and $\{\lambda^k_{\hat{i}}\} \to \infty$, we have
\begin{align}
u^Td_\bj&=\lim\limits_{k\to\infty} \sum\limits_{i\in\cI_=(\bx)}\lambda_i^k[\nabla \phi_i(\bx)+v^k_i-\nabla\psi_{i,j_i}(\bx)]^Td_\bj+(w^k)^Td_\bj\nn\\
&\le\limsup\limits_{k\to\infty} \sum\limits_{i\in\cI_=(\bx)}\lambda_i^k[\nabla \phi_i(\bx)^Td_\bj +\zeta'_i(\bx;d_\bj)-\nabla\psi_{i,j_i}(\bx)^Td_\bj]\nn\\
&\leq \limsup\limits_{k\to\infty} \lambda_{\hat{i}}^k[\nabla \phi_{\hat{i}}(\bx)^Td_\bj+\zeta'_{\hat{i}}(\bx;d_\bj)-\nabla\psi_{\hat{i},j_{\hat{i}}}(\bx)^Td_\bj]
=-\infty,\nn
\end{align}
which contradicts the fact that $u^Td_\bj$ is a constant. It follows that $P_\bj(\bx)$ is closed.

Secondly, we assume that $X$ is a polyhedral set, $\phi_{i}$ is affine and $\zeta_i$ is piecewise affine on $X$ for every $i\in\mathcal{I}_=(\bar{x})$. It follows that
$\partial \zeta_i(\bx) = \conv(\{u^{i,1}, \ldots, u^{i,m_i}\})$ for some vectors $u^{i,j}$ with $1\le j\le m_i$, where $\conv(\cdot)$ denotes the convex hull of the associated set. Let
\[
Q_{\bj}(\bx) =\left\{\sum\limits_{i\in\cI_=(\bx)}\lambda_i[\nabla \phi_i(\bx)+v_i-\nabla\psi_{i,j_i}(\bx)]: \lambda_i\geq0,  v_i\in\partial\zeta_i(\bx)\right\}.
\]
One can observe that
\[
Q_{\bj}(\bx) =
\left\{\sum\limits_{i\in\cI_=(\bx)}\left(\lambda_i[\nabla \phi_i(\bx)-\nabla\psi_{i,j_i}(\bx)]+\sum^{m_i}_{j=1} t_{ij} u^{i,j}\right): \sum^{m_i}_{j=1} t_{ij}  = \lambda_i,\lambda_i \geq 0,  t_{ij}  \ge 0
\right\}
\]
and hence it is a polyhedral cone. Since $X$ is a polyhedral set, ${\cal N}_X(\bx)$ is also
a polyhedral cone. Notice that $P_{\bj}(\bx)=Q_{\bj}(\bx)+{\cal N}_X(\bx)$. It then follows that $P_{\bj}(\bx)$ is a polyhedral cone and hence it is a closed convex cone.

Finally, since $P_{\bj}(\bx)$ is closed and $[P_\bj(\bx)]^\o=C_\bj(\bx)$,  we conclude from  \cite[Theorem 14.1]{Rockafellar} that $P_\bj(\bx)=([P_\bj(\bx)]^\o)^\o=[C_\bj(\bx)]^\o$.
\end{proof}

We are now ready to prove Theorem \ref{thm:KKT}.

\begin{proof}
($\Rightarrow$). Suppose that $\bx$ is a B-stationary point of problem \eqref{e1}. Let $j_0\in\cJ_0(\bx)$ and $\bj=(j_1,\cdots,j_I)\in\cJ(\bx)$ be arbitrarily chosen. By the above assumption on $\bx$, it follows from Theorem \ref{suff-cond-B} that $F'(\bx; d) \ge 0$ for all $d \in \cT_{Y_{\mathbb{j}}(\bar{x})}(\bar{x})$. It together with \eqref{Danskin} implies that for all $d \in \cT_{Y_{\mathbb{j}}(\bar{x})}(\bar{x})$,
\beq\label{p-first-ord}
\nabla\phi_0(\bx)^Td+\zeta'_0(\bx;d)-\nabla\psi_{0,j_0}(\bx)^Td
\geq\nabla\phi_0(\bx)^Td+\zeta'_0(\bx;d)-\max\limits_{j\in \mathcal{J}_0(x)}\nabla\psi_{0,j}(\bx)^Td=F'(\bx; d) \ge 0.
\eeq
For convenience, let $P_\bj=P_\bj(\bx)$ and $S_{j_0}=\{-\nabla\phi_0(\bx)-v_0+\nabla\psi_{0,j_0}(\bx): v_0\in\partial \zeta_0(\bx)\}$, where $P_\bj(\bx)$ is defined in \eqref{PJ-cone}. We next show that $P_\bj\cap S_{j_0}\neq\emptyset$. Suppose for contradiction that $P_\bj\cap S_{j_0}=\emptyset$. This, together
with the facts that $P_\bj$ and $S_{j_0}$ are nonempty closed convex sets and $S_{j_0}$  is bounded, implies that there exists some $\bd\in\Re^n$ such that $u^T\bd\le0$ for any $u\in P_\bj$ and $u^T\bd\ge1$ for any $u\in S_{j_0}$. Hence, one has
$\bd\in[P_\bj]^\o$.  By $\cI_=(\bx)\neq\emptyset$ and the above assumption on $\bx$, it follows from Propositions \ref{consqua} and \ref{polar-cone} that $[P_\bj]^\o=\cT_{Y_{\bj}(\bx)}(\bx)$, which along with $\bd\in[P_\bj]^\o$  implies that $\bd\in \cT_{Y_{\bj}(\bx)}(\bx)$. In addition, since $u^T\bd\ge1$ for any $u\in S_{j_0}$, one has
\[1\leq\min\limits_{u\in S_{j_0}} u^T\bd
=\min\limits_{v_0\in \partial\zeta_0(\bx)} \{-[\nabla\phi_0(\bx)+v_0-\nabla\psi_{0,j_0}(\bx)]^T\bd\}
=-[\nabla\phi_0(\bx)^T\bd+\zeta'_0(\bx;\bd)-\nabla\psi_{0,j_0}(\bx)^T\bd].
\]
Hence, we have that $\nabla\phi_0(\bx)^T\bd+\zeta'_0(\bx;\bd)-\nabla\psi_{0,j_0}(\bx)^T\bd <0$ and $\bd\in\cT_{Y_{\bj}(\bx)}(\bx)$, which contradict (\ref{p-first-ord}). We thus conclude that $P_\bj\cap S_{j_0}\neq\emptyset$. Using this relation, \eqref{PJ-cone}, the definitions of $S_{j_0}$ and $\cI_=(\bx)$, and letting $\lambda^{j_0,\bj}_i=0$ for all $i\not\in\cI_=(\bx)$, we easily see that \eqref{kkt-1} and \eqref{kkt-2} hold.

($\Leftarrow$). Suppose that for any $j_0\in\cJ_0(\bx)$ and $\bj=(j_1,\ldots,j_I)\in\cJ(\bx)$, there exists a vector of Lagrangian multipliers $\lambda^{j_0,\bj}=(\lambda_1^{j_0,\bj},\ldots,\lambda^{j_0,\bj}_I)$ such that \eqref{kkt-1} and \eqref{kkt-2} hold. In view of  \eqref{kkt-1}, one has $\lambda_i^{j_0,\bj} \ge 0$ for every $i\in\cI_=(\bx)$  and
$\lambda_i^{j_0,\bj}=0$ for any $i\not\in\cI_=(\bx)$. By \eqref{kkt-2}, there exist $v^{j_0,\bj}_i\in\partial\zeta_i(\bx)$ for all $0\le i \le I$ and $w^{j_0,\bj}\in\mathcal{N}_X(\bx)$ such that
\[\nabla\phi_0(\bx)+v^{j_0,\bj}_0-\nabla\psi_{0,j_0}(\bx)+\sum\limits_{i=1}^I\lambda_i^{j_0,\bj}[\nabla\phi_i(\bx)+v^{j_0,\bj}_i-\nabla\psi_{i,j_i}(\bx)]+w^{j_0,\bj}=0.\]
These, together with the facts that $\lambda_i^{j_0,\bj} \ge 0, \forall i\in\cI_=(\bx)$  and $\lambda_i^{j_0,\bj}=0,\forall i\not\in\cI_=(\bx)$, imply that for any $j_0\in\cJ_0(\bx)$, $\bj\in\cJ(\bx)$ and $d\in\Re^n$,
\begin{align}
\nabla\phi_0(\bx)^Td+\zeta_0'(\bx;d)-\nabla\psi_{0,j_0}(\bx)^Td&\geq\nabla\phi_0(\bx)^Td+(v^{j_0,\bj}_0)^Td-\nabla\psi_{0,j_0}(\bx)^Td\nn\\
&=-\sum\limits_{i=1}^I\lambda_i^{j_0,\bj}[\nabla\phi_i(\bx)+v^{j_0,\bj}_i-\nabla\psi_{i,j_i}(\bx)]^Td-(w^{j_0,\bj})^Td\nn\\
&\geq-\sum\limits_{i\in\cI_=(\bx)}\lambda_i^{j_0,\bj}[\nabla\phi_i(\bx)^Td+\zeta_i'(\bx;d)-\nabla\psi_{i,j_i}(\bx)^Td]-(w^{j_0,\bj})^Td.\label{p-first-ord-1}
\end{align}
By the above assumption on $\bx$, it follows from Proposition \ref{consqua} that $C_\bj(\bx)=\cT_{Y_\bj(\bx)}(\bx)$ for all $\bj\in\cJ(\bx)$. Hence,  one has that $d\in \cT_X(\bx)$ and $\nabla\phi_i(\bx)^Td+\zeta_i'(\bx;d)-\nabla\psi_{i,j_i}(\bx)^Td \le 0$ for any $d\in\cT_{Y_\bj(\bx)}(\bx)$. In view of these, \eqref{p-first-ord-1}, $w^{j_0,\bj}\in\mathcal{N}_X(\bx)$ and $\lambda_i^{j_0,\bj} \ge 0$ for every $i\in\cI_=(\bx)$, one has $
\nabla\phi_0(\bx)^Td+\zeta_0'(\bx;d)-\nabla\psi_{0,j_0}(\bx)^Td \ge 0$ for any $d\in\cT_{Y_\bj(\bx)}(\bx)$.
It then follows that for any $\mathbb{j}\in\mathcal{J}(\bar{x})$ and $d\in\cT_{Y_\bj(\bx)}(\bx)$,
\begin{align}
F'(\bx;d)&=\nabla\phi_0(\bx)^Td+\zeta_0'(\bx;d)-\max\limits_{j_0\in\cJ_0(\bx)}\nabla\psi_{0,j_0}(\bx)^Td\nn\\
&=\min\limits_{j_0\in\cJ_0(\bx)}\left\{\nabla\phi_0(\bx)^Td+\zeta_0'(\bx;d)-\nabla\psi_{0,j_0}(\bx)^Td\right\}\geq0.\nn
\end{align}
 Hence, we derive from Theorem \ref{suff-cond-B} that $\bx$ is a B-stationary point of problem \eqref{e1}.
\end{proof}

\end{document}